\documentclass{amsart}

\usepackage{amsthm}
\usepackage{amssymb}
\usepackage{amsmath}
\usepackage{amsfonts}
\usepackage{amsrefs}
\usepackage{textcomp}
\usepackage[T1]{fontenc}
\usepackage[utf8x]{inputenc}
\usepackage{textcomp}
\usepackage{wasysym}
\usepackage{stmaryrd}
\usepackage{esint}
\usepackage[all]{xy}
\usepackage{graphicx}
\usepackage{bbm}

\newtheorem{tw}{Theorem}[section]
\newtheorem{dfn}[tw]{Definition}
\newtheorem{uw}[tw]{Remark}
\newtheorem{prz}[tw]{Example}
\newtheorem{lem}[tw]{Lemma}
\newtheorem{stw}[tw]{Proposition}
\newtheorem{wn}[tw]{Corollary}
\newtheorem{zad}{Zadanie}
\newtheorem{hip}{Question}

\newtheorem*{ozn}{Notation}
\newtheorem*{dd}{Proof}
\newtheorem*{roz}{Rozwiązanie}
\newtheorem*{ak}{Acknowledgements}

\let\olddfn\dfn
\renewcommand{\dfn}{\olddfn\normalfont}

\let\oldozn\ozn
\renewcommand{\ozn}{\oldozn\normalfont}

\let\oldlem\lem
\renewcommand{\lem}{\oldlem\normalfont}

\let\oldstw\stw
\renewcommand{\stw}{\oldstw\normalfont}

\let\olduw\uw
\renewcommand{\uw}{\olduw\normalfont}

\let\oldwn\wn
\renewcommand{\wn}{\oldwn\normalfont}

\let\oldprz\prz
\renewcommand{\prz}{\oldprz\normalfont}

\let\olddd\dd
\renewcommand{\dd}{\olddd\normalfont}

\let\oldroz\roz
\renewcommand{\roz}{\oldroz\normalfont}

\let\oldak\ak
\renewcommand{\ak}{\oldak\normalfont}

\let\oldzad\zad
\renewcommand{\zad}{\oldzad\normalfont}

\let\oldhip\hip
\renewcommand{\hip}{\oldhip\normalfont}

\def\Xint#1{\mathchoice
   {\XXint\displaystyle\textstyle{#1}}%
   {\XXint\textstyle\scriptstyle{#1}}%
   {\XXint\scriptstyle\scriptscriptstyle{#1}}%
   {\XXint\scriptscriptstyle\scriptscriptstyle{#1}}%
   \!\int}
\def\XXint#1#2#3{{\setbox0=\hbox{$#1{#2#3}{\int}$}
     \vcenter{\hbox{$#2#3$}}\kern-.5\wd0}}

\def\dashint{\Xint-}

\usepackage{geometry}
\author{Wojciech G\'{o}rny}

\address{W. G\'{o}rny: Faculty of Mathematics, Informatics and Mechanics, University of Warsaw, Warsaw, Poland.}

\email{w.gorny@mimuw.edu.pl}

\date{\today}


\title{Least gradient problem with respect to a non-strictly convex norm}


\begin{document}

\begin{abstract}
We study the planar least gradient problem with respect to an anisotropic norm $\phi$ for continuous boundary data. We prove existence of minimizers for strictly convex domains $\Omega$. Furthermore, we inspect the issue of uniqueness and regularity of minimizers only in terms of the modes of convexity of $\phi$ and $\Omega$. The results are independent from the regularity of $\phi$.
\end{abstract}

\maketitle

\section{Introduction}

We are interested in the issue of existence and uniqueness of minimizers to the least gradient problem in the following setting:

\begin{equation}\label{problem}\tag{ALGP}
\min \{ \int_\Omega \phi(Du), \quad u \in BV(\Omega), \quad u|_{\partial\Omega} = f \},
\end{equation}
where we assume that the metric integrand $\phi$ is convex without any additional regularity assumptions. Note that $\phi$ depends only on the direction of the derivative. As for discontinuous boundary data there can be multiple solutions or no solutions at all even in the isotropic case, see \cite{MRL} and \cite{ST} respectively, we assume continuity of the boundary data $f \in C(\partial\Omega)$. Moreover, throughout this paper we assume that $\Omega \subset \mathbb{R}^2$.

In the classical least gradient problem, when $\phi$ is the Euclidean norm, existence, regularity, and uniqueness of minimizers depends on the geometry of the set $\Omega$. Here the situation is slightly more complicated, as we have additionally the interplay between the shapes of $\Omega$ and the unit ball in the anisotropic norm $B_\phi(0,1)$; our goal is to explore this relationship. We divide our reasoning into two stages:

(1) Suppose that the unit ball $B_{\phi}(0,1)$ is strictly convex. Then, regardless of the regularity of $\phi$, we are able to prove existence and uniqueness of minimizers for strictly convex $\Omega$ and obtain regularity estimates in terms of the modulus of continuity of the boundary data, which are independent on the choice of $\phi$.

(2) Suppose that the unit ball $B_{\phi}(0,1)$ has flat facets. Then, under stronger assumptions on $\Omega$, we use the regularity estimates from the strictly convex case to prove existence of a single minimizer with the same regularity. However, we lose uniqueness of minimizers and the additional minimizers may have regularity no better than $BV(\Omega) \cap L^\infty(\Omega)$.

Let us stress that in the anisotropic least gradient problem known results on uniqueness of minimizers for continuous boundary data depend not only on the geometry of $\Omega$, but also on the regularity of $\phi$. For instance, the uniqueness proof in \cite{JMN} is based on a maximum principle and requires uniform convexity and a condition slightly weaker than $W^{3,\infty}$ regularity of an elliptic metric integrand away from $\{ \xi = 0 \}$; for a precise assumption, see \cite[Theorem 1.2]{JMN}. In the course of this paper, we are going to relax the assumptions on the regularity of $\phi$ in the case when $\phi$ depends only on the second (directional) variable in order to be able to deal with non-strictly convex metric integrands.

This paper is organized as follows: in Section 2 we recall the basic definitions and facts concerning anisotropic $BV$ spaces and $\phi-$least gradient functions. In Section 3 we investigate the functional
$$ F_{\phi}(v) = \int_\Omega |Dv|_{\phi} + \int_{\partial\Omega} \phi(x,\nu^\Omega) |Tv - f| d\mathcal{H}^{N-1}$$
and (without the restrictions as to the form of $\phi$ or the dimension) prove that if metric integrands $\phi_n \rightarrow \phi$ uniformly in $C(\overline{\Omega} \times \partial B(0,1))$, then the functionals $F_{\phi_n}$ $\Gamma-$converge to $F_\phi$. In particular, it provides a stability result in the spirit of Miranda's theorem, \cite{Mir}, which states that a sequence of (isotropic) least gradient functions converges to a least gradient function. 

In Section 4 we explore the results of Jerrard, Nachman and Tamasan, \cite{JMN}, and extend the framework under which they are valid. In \cite{JMN} the two most important assumptions are the {\it barrier condition}, see Definition \ref{def:barrier}, which is essential in existence proofs, and uniform convexity and quite strong regularity of the metric integrand $\phi$, which is used in uniqueness proofs. Here, we prove existence and uniqueness of minimizers for strictly convex $B_\phi(0,1)$ regardless of the regularity of $\phi$; see Theorem \ref{tw:strict}. Moreover, we prove a regularity estimate for the minimizer depending only on $\Omega$ and regularity of the boundary data, see Proposition \ref{prop:modulusofcontinuity}; it does not depend on the regularity of $\phi$.

In Section 5 we show that only some results can be extended to the case when $B_\phi(0,1)$ is not strictly convex. In Theorem \ref{tw:existence} we prove existence of a minimizer, which has the same regularity as if $B_\phi(0,1)$ was strictly convex, provided that $\Omega$ is uniformly convex; in Theorem \ref{tw:existence2} we prove existence of a minimizer for strictly convex $\Omega$. However, we lose uniqueness of minimizers and not all minimizers reflect the same regularity and we also include examples (see Proposition \ref{stw:example}) such that the solution $u$ has regularity no better than $BV(\Omega)$ even if the boundary data are smooth.

Finally, in Section 6 we discuss which results in the isotropic case can be extended to the anisotropic case due to the newly estabilished results. Moreover, in the various stages of the reasoning, we interrupt it to focus solely on the geometry underneath these results: we discuss where stronger modes of convexity come into play and which sets $\Omega$ satisfy the barrier condition.

\section{Preliminaries}

This Section is divided into two principal parts: firstly, we recall some general facts about $BV$ spaces with respect to an anisotropic norm $\phi$. Our main point of reference for this subsection is \cite{AB}. Then we define functions of $\phi-$least gradient and recall some of their properties.

\subsection{$BV$ spaces with respect to an anisotropic norm}

\begin{dfn}
Let $\Omega \subset \mathbb{R}^N$ be an open bounded set with Lipschitz boundary. A continuous function $\phi: \overline{\Omega} \times \mathbb{R}^N \rightarrow [0, \infty)$ is called a metric integrand, if it satisfies the following conditions: \\
\\
$(1)$ $\phi$ is convex with respect to the second variable for a.e. $x \in \overline{\Omega}$; \\
$(2)$ $\phi$ is homogeneous with respect to the second variable, i.e.
\begin{equation*}
\forall \, x \in \overline{\Omega}, \quad \forall \, \xi \in \mathbb{R}^N, \quad \forall \, t \in \mathbb{R} \quad \phi(x, t \xi) = |t| \phi(x, \xi);
\end{equation*}
$(3)$ $\phi$ is bounded in $\overline{\Omega}$, i.e.
\begin{equation*}
\exists \, \Gamma > 0 \quad \forall \, x \in \overline{\Omega}, \quad \forall \, \xi \in \mathbb{R}^N \quad 0 \leq \phi(x, \xi) \leq \Gamma |\xi|.
\end{equation*}
\\
The above conditions are not enough to define the anisotropic total variation in a way that recovers some properties of the classical total variation, so we will additionally assume that \\
\\
$(4)$ $\phi$ is elliptic in $\overline{\Omega}$, i.e.
\begin{equation*}
\exists \, \lambda > 0 \quad \forall \, x \in \overline{\Omega}, \quad \forall \, \xi \in \mathbb{R}^N \quad \lambda |\xi| \leq \phi(x, \xi).
\end{equation*}
\end{dfn}

These conditions apply to most cases considered in the literature, such as the classical least gradient problem, i.e. $\phi(x, \xi) = |\xi|$ (see \cite{SWZ}), the weighted least gradient problem, i.e. $\phi(x, \xi) = g(x) |\xi|$ (see \cite{JMN}), where $g \geq c > 0$, and $l_p$ norms for $p \in [1, \infty]$, i.e. $\phi(x, \xi) = \| \xi \|_p$ (see \cite{Gor1}).

\begin{dfn}
The polar function of $\phi$ is $\phi^0: \overline{\Omega} \times \mathbb{R}^N \rightarrow [0, \infty)$ defined as

\begin{equation*}
\phi^0 (x, \xi^*) = \sup \, \{ \langle \xi^*, \xi \rangle : \, \xi \in \mathbb{R}^N, \, \phi(x, \xi) \leq 1 \}.
\end{equation*} 
\end{dfn}

\begin{dfn}
Let $\phi$ be a metric integrand continuous and elliptic in $\overline{\Omega}$. For a given function $u \in L^1(\Omega)$ we define its $\phi-$total variation in $\Omega$ by the formula:

\begin{equation*}
\int_\Omega |Du|_\phi = \sup \, \{ \int_\Omega u \, \mathrm{div} \, \mathbf{z} \, dx : \, \phi^0(x,\mathbf{z}(x)) \leq 1 \, \, \, \text{a.e.}, \quad \mathbf{z} \in C_c^1(\Omega)  \}.
\end{equation*}
In the literature the total variation is sometimes instead denoted by $\int_\Omega \phi(x, Du)$. We say that $u \in BV_\phi(\Omega)$ if its $\phi-$total variation is finite. Similarly, we define the $\phi-$perimeter of a set $E$ to be

$$P_\phi(E, \Omega) = \int_{\Omega} |D\chi_E|_\phi.$$
If $P_\phi(E, \Omega) < \infty$, we say that $E$ is a set of bounded $\phi-$perimeter in $\Omega$.
\end{dfn}

\begin{uw}
By properties $(3)$ and $(4)$ of a metric integrand we have that $\lambda \int_\Omega |Du| \leq \int_\Omega |Du|_\phi \leq \Gamma \int_\Omega |Du|$. In particular, $BV_\phi(\Omega) = BV(\Omega)$ as sets; however, they are equipped with different (but equivalent) norms and corresponding strict topologies.
\end{uw}

The spaces $BV_{\phi}(\Omega)$ defined as above satisfy the same basic properties as the isotropic space $BV(\Omega)$: we recover lower semicontinuity of the $\phi-$total variation with respect to $L^1$ convergence, the isoperimetric inequality and the co-area formula. Moreover, we recover the approximation by smooth functions in the strict topology:

\begin{uw}\label{lem:scisaniz}
Suppose that $\Omega$ is an open bounded set with Lipschitz boundary and $\phi$ is a metric integrand continuous and elliptic in $\overline{\Omega}$. Let $v \in BV_\phi(\Omega)$ and $Tv = f$. Then there exists a sequence $v_n \in C^\infty(\Omega) \cap BV(\Omega)$ such that $v_n \rightarrow v$ strictly in $BV_\phi(\Omega)$ and $Tv_n = f$ (in the isotropic case, see \cite[Corollaries 1.17, 2.10]{Giu}).

The approximation by smooth functions in the strict topology entails that we may approximate sets of bounded $\phi-$perimeter both in the Lebesgue measure and in $\phi-$perimeter by open sets with smooth boundary (also with respect to some given boundary conditions). For the proof in the isotropic case, see \cite[Theorem 3.42]{AFP}. \qed
\end{uw}

Throughout most of this paper, we will use the following integral representation of the $\phi-$total variation (\cite{AB}, \cite{JMN}):

\begin{stw}\label{stw:repcalkowa}
Let $\varphi: \overline{\Omega} \times \mathbb{R}^N \rightarrow \mathbb{R}$ be a metric integrand. Then we have an integral representation:

\begin{equation*}
\int_\Omega |Du|_\phi = \int_\Omega \phi(x, \nu^u(x)) \, |Du|,
\end{equation*}
where $\nu^u$ is the Radon-Nikodym derivative $\nu^u = \frac{d Du}{d |Du|}$. If we take $u$ to be a characteristic function of a set $E$ with a $C^1$ boundary, we have

\begin{equation*}
P_\phi(E, \Omega) = \int_\Omega \phi(x, \nu_E) \, d \mathcal{H}^{N-1},
\end{equation*}
where $\nu(x)$ is the (Euclidean) unit vector normal to $\partial E$ at $x \in \partial E$. Let us also denote by $\tau(x)$ the unit vector tangent to $\partial E$ at $x \in \partial E$. \qed
\end{stw}

\subsection{$\phi-$least gradient functions}

Now, we turn our attention to the precise formulation of Problem (\ref{problem}). Then we recall several known properties of the minimizers and a few results concerning the isotropic case, i.e. $\phi(x, \xi) = |\xi|$.

\begin{dfn}
Let $\Omega \subset \mathbb{R}^N$ be an open bounded set with Lipschitz boundary. We say that $u \in BV_\phi(\Omega)$ is a function of $\phi-$least gradient, if for every compactly supported $v \in BV_\phi(\Omega)$ we have

\begin{equation*}
\int_\Omega |Du|_\phi \leq \int_\Omega |D(u + v)|_\phi.
\end{equation*}
If $\phi$ is a metric integrand with continuous extension to $\mathbb{R}^N$, we may instead assume that $v$ is a $BV_\phi$ function with zero trace on $\partial \Omega$; see \cite[Proposition 3.16]{Maz}. We say that $u$ is a solution to Problem (\ref{problem}), the anisotropic least gradient problem with boundary data $f$, if $u$ is a function of $\phi-$least gradient and $Tu = f$.
\end{dfn}

Both in the isotropic and anisotropic case existence and uniqueness of minimizers depend on the geometry of $\Omega$. Suppose that the boundary data are continuous. In the isotropic case, the necessary and sufficient condition was introduced in \cite{SWZ} and in two dimensions it is equivalent to strict convexity of $\Omega$. In the anisotropic case, a sufficient condition (see \cite[Theorem 1.1]{JMN}) is the {\it barrier condition}:

\begin{dfn}\label{def:barrier}
(\cite[Definition 3]{JMN}) Let $\Omega \subset \mathbb{R}^N$ be an open bounded set with Lipschitz boundary. Suppose that $\phi$ is an elliptic metric integrand. We say that $\Omega$ satisfies the barrier condition if for every $x_0 \in \partial\Omega$ and sufficiently small $\varepsilon > 0$, if $V$ minimizes $P_\phi(\cdotp; \mathbb{R}^N)$ in 

$$ \{ W \subset \Omega: W \backslash B(x_0, \varepsilon) = \Omega \backslash B(x_0, \varepsilon) \} $$
then
$$\partial V \cap \partial\Omega \cap B(x_0, \varepsilon) = \emptyset.$$
In the isotropic case $\phi(x, \xi) = \| \xi \|_2$ this is equivalent, at least for sets with $C^2$ boundary, to the condition introduced in \cite{SWZ}.
\end{dfn}

Before we proceed, we need one additional result that relates functions of $\phi-$least gradient and $\phi-$minimal sets, i.e. sets, the characteristic functions of which are of $\phi-$least gradient. It may be seen as an anisotropic version of the classical result by Bombieri, de Giorgi and Giusti, see \cite{BGG}. Its proof in both directions is based on the the co-area formula.

\begin{stw}\label{stw:anizobgg}
(\cite[Theorem 3.19]{Maz}) Let $\Omega \subset \mathbb{R}^N$ be an open bounded set with Lipschitz boundary. Assume that the metric integrand $\phi$ has a continuous extension to $\mathbb{R}^N$. Take $u \in BV_\phi(\Omega)$. Then $u$ is a function of $\phi-$least gradient in $\Omega$ if and only if $\chi_{\{ u > t \}}$ is a function of $\phi-$least gradient for almost all $t \in \mathbb{R}$. \qed
\end{stw}

Finally, as $\phi-$least gradient functions are $BV$ functions, they are defined up to a set of measure zero, we have to choose a proper representative if we want to state any regularity results. In this paper, following \cite{SWZ}, we employ the convention that a set of a bounded perimeter consists of all its points of positive density.

\section{$\Gamma-$convergence}

We start with recalling the notion of $\Gamma-$convergence:

\begin{dfn}
Let $F, F_n: X \rightarrow [0, \infty]$ be a sequence of functionals on a topological space $X$. We say that the sequence $F_n$ $\Gamma-$converges to $F$, what we denote by $\Gamma-\lim_{n \rightarrow \infty} F_n = F$, if the following two conditions are satisfied: \\
(1) For every sequence $x_n \in X$ such that $x_n \rightarrow x$ in $X$ we have 
$$F(x) \leq \liminf_{n \rightarrow \infty} F_n(x_n);$$
(2) For every $x \in X$ there exists a sequence $x_n \rightarrow x$ in $X$ such that
$$F(x) \geq \limsup_{n \rightarrow \infty} F_n(x_n).$$
An important property of $\Gamma-$convergence is that cluster points of minimizers of $F_n$ are minimizers of $F$.
\end{dfn}

We define the following functional (which is a relaxation of the total variation functional with respect to Dirichlet boundary data, see \cite{Maz}): 

$$ F_{\phi}(v) = \int_\Omega |Dv|_{\phi} + \int_{\partial\Omega} \phi(x,\nu^\Omega) |Tv - f| d\mathcal{H}^{N-1}.$$

Now we state the main result in this Section. The main idea behind it is extending Miranda's theorem, see \cite{Mir}, which states that a sequence of least gradient functions convergent in $L^1(\Omega)$ converges to a least gradient function. The Theorem below allows us to consider $\phi_n$-least gradient functions, where the anisotropic norm $\phi_n$ is not fixed and changes with $n$. This enables us to prove existence results in the anisotropic least gradient problem, when $\partial B_\phi(0,1)$ is not strictly convex, as we need to approximate the anisotropic norm using strictly convex metric integrands.

\begin{tw}\label{tw:gammaconvergence}
Let $\phi$ and $\phi_n$ be metric integrands such that $\phi_n \rightarrow \phi$ in $C(\overline{\Omega} \times \partial B(0,1))$. Then the sequence of functionals $F_{\phi_n}$ $\Gamma-$converges (with respect to the $L^1$ convergence) to the functional $F_\phi$.
\end{tw}

\begin{dd}
(1) We show that for any sequence $u_n \rightarrow u$ in $L^1(\Omega)$ we have $F_\phi(u) \leq \liminf_{n \rightarrow \infty} F_{\phi_n}(u_n)$.

Denote by $F_{l_2}$ the functional $F_\psi$, where $\psi$ is the isotropic norm. If $\liminf_{n \rightarrow \infty} F_{\phi_n}(u_n) = \infty$, then the inequality is obvious. Assume that this number is finite and take the subsequence, still denoted by $u_{n}$, such that this limit is achieved. In particular, the (new) sequence $(u_{n})$ is bounded in $BV(\Omega)$ and thus $F_{l^2}(u_n)$ is bounded. Assume for now that $F_\phi(u) < \infty$. Then

$$ \liminf_{n \rightarrow \infty} (F_{\phi_n}(u_n) - F_\phi(u)) \geq \liminf_{n \rightarrow \infty} (F_{\phi_n}(u_n) - F_\phi(u_n)) + \liminf_{n \rightarrow \infty} (F_{\phi}(u_n) - F_\phi(u))$$
and by the lower semicontinuity of $F_\phi$ the second summand is nonnegative. Hence

$$ \liminf_{n \rightarrow \infty} F_{\phi_n}(u_n) \geq \liminf_{n \rightarrow \infty} (F_{\phi_n}(u_n) - F_\phi(u_n)) = 0,$$
because we have

$$ |F_{\phi_n}(u_n) - F_\phi(u_n)| = \int_\Omega (\phi_n(x,\nu^{u_n}) - \phi(x,\nu^{u_n}))|Du_n| + \int_{\partial\Omega} (\phi_n(x,\nu^\Omega) - \phi(x,\nu^\Omega)) |Tu_n - f| d\mathcal{H}^{N-1} \leq$$
$$ \leq \sup_{\overline{\Omega} \times \partial B(0,1)} |\phi_n - \phi| (\int_\Omega |Du_n| + \int_{\partial \Omega} |Tu_n - f|) = \sup_{\partial B(0,1)} |\phi_n - \phi| \, F_{l^2}(u_n) \leq \sup_{\partial B(0,1)} |\phi_n - \phi| \, M \rightarrow 0.$$

If $F_\phi(u) = \infty$, then $u \notin BV(\Omega)$. By the lower semicontinuity of the total variation for every approximating sequence $u_n \rightarrow u$ we also have $\liminf_{n \rightarrow \infty} |Du_n| = \infty$. But as the sequence $\phi_n$ converges uniformly, $\phi_n-$norms are uniformly equivalent to the Euclidean norm, so

$$ \liminf_{n \rightarrow \infty} F_{\phi_n}(u_n) \geq \liminf_{n \rightarrow \infty} \int_\Omega |Du_n|_{\phi_n} \geq \liminf_{n \rightarrow \infty} C \int_\Omega |Du_n| = \infty.$$

(2) We show that for any function $u \in L^1(\Omega)$ there exists a sequence $u_n \rightarrow u$ such that $F_\phi(u) \geq \limsup_{n \rightarrow \infty} F_{\phi_n}(u_n)$. In fact, it is enough to consider the constant sequence $u_n  = u$.

If $u \notin BV(\Omega)$, the inequality is obvious. If $u \in BV(\Omega)$, then all the integrals in the definitions of $F_{\phi}, F_{\phi_n}, F_{l^2}$ are convergent and we have

$$ |F_{\phi_n}(u) - F_\phi(u)| = \int_\Omega (\phi_n(x,\nu^u) - \phi(x,\nu^u))|Du| + \int_{\partial\Omega} (\phi_n(x,\nu^\Omega) - \phi(x,\nu^\Omega)) |Tu - f| d\mathcal{H}^{N-1} \leq$$
$$ \leq \sup_{\partial B(0,1)} |\phi_n - \phi| (\int_\Omega |Du| + \int_{\partial \Omega} |Tu - f|) = \sup_{\partial B(0,1)} |\phi_n - \phi| \, F_{l^2}(u) \rightarrow 0.$$
\qed
\end{dd}

The assumption that $\phi_n \rightarrow \phi$ in $C(\overline{\Omega} \times \partial B(0,1))$ is quite natural in this context: as metric integrands are $1-$homogenous in the second variable, it is sufficient to check convergence only on the unit sphere. Furthermore, as the following Example shows, we may not relax the assumption concerning uniform convergence in $\Omega$; however, as we can see in Proposition \ref{stw:pointuniform}, some form of uniform convergence in the second variable is guaranteed.

\begin{prz}
Let $\Omega = [-1,1]$, $f(-1) = 0$ and $f(1) = 1$. Let $\phi_n(x, p) = a_n(x) ||p||_{l^2}$, where $a_n \in C^{1,1}(\Omega)$ such that $a_n(x) \in [\frac{1}{2},1]$, $\min a_n = a_n(x_n) = \frac{1}{2}$, where $x_n \rightarrow 0$ and $a_n(x) = 1$ for $x \in [-1, x_n - \frac{1}{n}] \cup [x_n + \frac{1}{n},1]$. We note that $a_n \rightarrow a$ pointwise and in every $L^q([-1,1])$, $q < \infty$, where $a \equiv 1$; thus $\phi_n \rightarrow \phi$ pointwise and in $L^q(\overline{\Omega} \times \partial B(0,1))$, where $\phi(x,p) = \| p \|_{l^2}$.

Let $u_n = \chi_{[x_n,1]}$. We have that $u_n \rightarrow u = \chi_{[0,1]}$ in $L^1([-1,1])$. We notice that both $u_n$ and $u$ have trace $f$. Thus

$$ F_\phi(u) = \int_{[-1,1]} a(x)|Du|= \int_{[-1,1]} |Du| = 1 > \frac{1}{2} = a_n(x_n) = \int_{[-1,1]} a_n(x) |Du_n| = F_{\phi_n}(u_n).$$
In particular, the first condition in the definition of $\Gamma-$convergence is not satisfied. Thus uniform convergence in $\overline{\Omega}$ is a necessary condition for $\Gamma-$convergence of $F_{\phi_n}$ to $F_{\phi}$.
\end{prz}

\begin{stw}\label{stw:pointuniform}
Let $\phi_n, \phi$ be metric integrands such that $\phi_n(x, p) = \phi_n(p)$ and suppose that $\phi_n \rightarrow \phi$ pointwise. Then $\phi_n \rightarrow \phi$ in $C(\partial B(0,1))$.
\end{stw}

\begin{dd}
In a finite dimensional space any norm is equivalent to the isotropic norm, so each $\phi_n$ is continuous and attains its supremum on $\partial B(0,1)$. As $\phi_n$ is convex and $1-$homogenous, we have

$$ |\phi_n(x) - \phi_n(y)| \leq \phi_n(x-y) = \phi_n(\frac{x-y}{|x-y|}) |x-y| \leq (\sup_{\partial B(0,1)} \phi_n) |x-y|,$$

so $\phi_n$ is Lipschitz continuous with constant $\sup_{\partial B(0,1)} \phi_n$. We have two possibilities: \\

1. $\sup_n (\sup_{\partial B(0,1)} \phi_n) \leq M$. Take any subsequence $\phi_{n_k}$. By Arzela-Ascoli theorem it has a convergent subsequence $\phi_{n_{k_l}}$; by our assumption we have that $\phi_{n_{k_l}} \rightarrow \phi$ uniformly. As $C(\partial B(0,1))$ is a metric space, we have that $\phi_n \rightarrow \phi$ in $C(\partial B(0,1))$. \\

2. $\sup_n (\sup_{\partial B(0,1)} \phi_n) = + \infty$. Let $x_n$ be the point where $\phi_n$ attains its supremum on $\partial B(0,1)$. Take $\{ q_k \}$ to be the $\varepsilon-$net on $\partial B(0,1)$ for $\varepsilon = \frac{1}{2}$. Fix $n$ and let $q_k$ be such that $|q_k - x_n| \leq \frac{1}{2}$. By the Lipschitz continuity of $\phi_n$

$$ \phi_n(q_k) \geq \phi_n(x_n) - (\sup_{\partial B(0,1)} \phi_n)|q_k - x_n| \geq \frac{1}{2} \phi_n(x_n).$$
As the set $\{ q_k \}$ is finite, for some $k$ the sequence $\phi_n(q_k)$ is unbounded. But this is impossible, as $\phi_n(q_k) \rightarrow \phi(q_k)$. \qed
\end{dd}

\section{Strictly convex unit ball $B_\phi(0,1)$}

From now on, we introduce the following notation:

\begin{ozn}
Let $\Omega \subset \mathbb{R}^2$ be an open bounded set with Lipschitz boundary. We denote by $\phi$ (or $\phi_n$) a metric integrand on $\overline{\Omega} \times \mathbb{R}^2$ which depends only on the second variable, i.e. $\phi(x,\xi) = \phi(\xi)$. Moreover, we will denote by $l_p$ the metric integrand defined by the formula $\phi(x,\xi) = \| \xi \|_p$.
\end{ozn}

When necessary, we will additionally assume some form of convexity of $\Omega$. Our reasoning is divided into two main parts: in this Section, we explore the case when the unit ball $B_\phi(0,1)$ is strictly convex and in the next Section we explore the case when $B_\phi(0,1)$ has flat parts of the boundary.

This Section is divided into two subsections. Firstly, we explore how do the $\phi-$minimal sets look like and use this knowledge to infer that for continuous boundary data the minimizers exist and are unique if $B_\phi(0,1)$ is strictly convex. In the second part, we will prove that minimizers inherit some of the regularity of the boundary data. As in the isotropic least gradient problem (see \cite{SWZ}), existence and regularity of minimizers require respectively strict convexity and uniform convexity of $\Omega$; however, let us underline the fact that these results will not depend on the regularity of $\phi$.

\subsection{Existence and uniqueness of minimizers}

The next result states that line segments always are $\phi-$minimal surfaces regardless of the regularity of $\phi$. However, it does not imply that there are no other $\phi-$minimal surfaces. We also note that the proof heavily relies on the fact that $\Omega$ is two-dimensional.

\begin{stw}\label{prop:linesegment}
Suppose that $\Omega$ is convex. Let $p_1, p_2 \in \partial\Omega$. Then the set $E \subset \Omega$ such that $\partial E$ equals the line segment $\overline{p_1 p_2}$ is a $\phi-$minimal set.
\end{stw}

\begin{dd}
Firstly, we recall the Jensen inequality: if $\phi: X \rightarrow \mathbb{R}$ is a convex and lower semicontinuous function, $X$ is a separable Banach space and $\mu$ is a probability measure, then

$$ \phi(\int_X x d\mu) \leq \int_X \phi(x) d\mu.$$

Now, let $E$ be as above and $F$ be a set with $C^1$ boundary such that $\partial F$ is a curve with ends $p_1, p_2$. Firstly, we notice that by the fundamental theorem of calculus we have (assume that $\partial E$ and $\partial F$ are oriented from $p_1$ to $p_2$):

\begin{equation*}
\int_{\partial F} \tau(x) d\mathcal{H}^{1} = p_2 - p_1 = \mathcal{H}^{1}(\partial E) \tau_0,
\end{equation*}
where $\tau_0$ is the unit vector along the line segment $\overline{p_1 p_2}$. Thus, by rotating the above equation by $\frac{\pi}{2}$, we obtain

\begin{equation}\label{eq:nu0}
\int_{\partial F} \nu(x) d\mathcal{H}^{1} = \mathcal{H}^{1}(\partial E) \nu_0,  
\end{equation}
where $\nu_0$ is the unit vector normal to the line segment $\overline{p_1 p_2}$. Now, we calculate

$$ P_\phi(E, \Omega) = \int_{\partial E} \phi(\nu_0) \, d\mathcal{H}^{1} = \mathcal{H}^{1}(\partial E) \, \phi(\nu_0) = \mathcal{H}^{1}(\partial E) \phi(\frac{1}{\mathcal{H}^{1}(\partial E)} \int_{\partial F} \nu(x) \, d\mathcal{H}^{1}) = $$
$$ = \phi(\int_{\partial F} \nu(x) \, d\mathcal{H}^{1}) = \mathcal{H}^{1}(\partial F) \phi(\int_{\partial F} \nu(x) \, d \frac{\mathcal{H}^{1}}{\mathcal{H}^{1}(\partial F)}) \leq $$
$$ \leq \mathcal{H}^{1}(\partial F) \int_{\partial F} \phi(\nu(x)) \, d \frac{\mathcal{H}^{1}}{\mathcal{H}^{1}(\partial F)} = \int_{\partial F} \phi(\nu(x)) \, d \mathcal{H}^{1} = P_\phi(F, \Omega).$$
The first and the second equality come from the integral representation of $P_\phi(E, \Omega)$ and the fact that $\partial E$ is an line segment. The third equality follows by the formula on $\nu_0$ given in (\ref{eq:nu0}). The fourth and the fifth equality hold because $\phi$ is $1-$homogenous. The inequality follows from the Jensen inequality for the convex function $\phi$ and the probability measure $\frac{\mathcal{H}^{1}}{\mathcal{H}^{1}(\partial F)}$ on $\partial F$. The sixth equality again holds because $\phi$ is $1-$homogenous and the final equality is the integral representation of $P_\phi(F, \Omega)$. 

To end the proof, we notice that we can approximate any set $F$ in the strict topology (with respect to $\phi$) with smooth sets $\widetilde{F_n}$ with the same boundary values. The boundary of $\widetilde{F_n}$ has the following structure:

$$ \partial \widetilde{F_n} = \, \bigcup_j \, \Gamma_n^j \, \cup \, \bigcup_i \, \Gamma_{i,n},$$
where $\Gamma_n^j$ are smooth curves from $p_1$ to $p_2$ and $\Gamma_{i,n}$ are closed loops. By Jordan curve theorem the union over $j$ is nonempty, as some part of the boundary has to disconnect the set with trace values equal to $0$ and $1$. Thus we may introduce the set $F_n$, whose boundary is $\Gamma_n^1$ with trace values the same as $E$. Then, by the previous calculation, we have

$$ P_\phi(E, \Omega) \leq P_\phi(F_n, \Omega) \leq P_\phi(\widetilde{F_n}, \Omega) \rightarrow P_\phi(F, \Omega),$$
so $E$ is a $\phi$-minimal set. \qed
\end{dd}

Let us stress that the above result does not state that the line segment is the only anisotropic minimal surface, only that it minimizes the perimeter for some special boundary conditions. Our goal is to show that the boundary of any $\phi-$minimal set $E$ is an at most countable union of line segments if $B_\phi(0,1)$ is strictly convex; this is not true for non-strictly convex unit ball and this case will be discussed separately in Section 5.

We proceed in two steps: firstly, we show that if $B_\phi(0,1)$ is strictly convex, then in the setting of Proposition \ref{prop:linesegment} the set $E$ is the only $\phi-$minimal set with respect to its boundary data. The simplest approach, which is to look closely at the proof of Proposition \ref{prop:linesegment} and check if the Jensen inequality is strict, fails if the boundary of a minimal set is not $C^1$. We will go around this problem by taking $F$, another $\phi-$minimal set relative to the same boundary data, approximate $F$ in the strict topology with smooth sets $F_n$ and prove a common lower bound for $P(F_n, \Omega)$. This approach is formalized in the following Proposition.

\begin{stw}\label{stw:strict}
Suppose that $\Omega$ is convex and suppose that $B_\phi(0,1)$ is strictly convex. Let $p_1, p_2 \in \partial\Omega$. Let $E \subset \Omega$ be a set such that $\partial E$ equals the line segment $\overline{p_1 p_2}$. Then $E$ is $\phi$-minimal. Moreover, if $F \subset \Omega$ is $\phi$-minimal and $T\chi_F = T\chi_E$, then $E = F$.
\end{stw}

\begin{dd}
Suppose that $F$ is another $\phi-$minimal set with the same trace and $l^2(E \Delta F) > 0$. By \cite[Corollary 2.6]{Gor2} the sets $E \cup F$ and $E \cap F$ are also $\phi-$minimal sets with the same trace and at least one of them differs from $E$ on a set of positive measure; thus we may replace $F$ by this set and without loss of generality assume that $E \subset F$.

We approximate $F$ in the strict topology (with respect to $\phi$) by smooth sets $\widetilde{F_n}$ with the same trace. In particular, we have a decomposition

$$ \partial \widetilde{F_n} = \Gamma_n \cup \, \bigcup_{i} \, \Gamma_{i, n},$$
where $\Gamma_n$ is a smooth curve from $p_1$ to $p_2$ and $\Gamma_{i,n}$ are closed loops. We notice that since it is an approximation of a minimal set, for sufficiently large $n$ in this decomposition there is only one smooth curve from $p_1$ to $p_2$, because by Proposition \ref{prop:linesegment} each of such curves has anisotropic length bounded from below by $P_\phi(E, \Omega)$. 

Now, let $A_{i,n}$ be sets enclosed by the loops $\Gamma_{i,n}$ and take $A_n = \bigcup_i A_{i,n}$. We want to modify the sets $\widetilde{F_n}$ in order to eliminate the closed loops: let $F_n$ be a set with the same trace as $F$ such that its boundary is $\Gamma_n$. We observe that $\widetilde{F_n} \Delta F_n \subset A_n$.

Since $\widetilde{F_n}$ is an approximation of $F$ in the strict topology and the set $F$ is $\phi$-minimal, we have

$$ P_\phi(F, \Omega) \leq P_\phi(F_n, \Omega) \leq P_\phi(F_n, \Omega) + P_\phi(A_n, \Omega) \leq P_\phi(\widetilde{F_n}, \Omega) \rightarrow P_\phi(F, \Omega).$$
Thus $P_\phi(F_n, \Omega) \rightarrow P_\phi(F, \Omega)$ and $P_\phi(A_n, \Omega) \rightarrow 0$; by the isoperimetric inequality also $l^2(A_n) \rightarrow 0$. In particular $l^2(\widetilde{F_n} \Delta F_n) \rightarrow 0$, so also 

$$l^2(F_n \Delta F) \leq l^2(F_n \Delta \widetilde{F_n}) + l^2(\widetilde{F_n} \Delta F) \rightarrow 0.$$
Thus $F_n$ is another approximation of $F$ in the strict topology with respect to $\phi$. 

Now, we are going to prove a lower bound on the anisotropic perimeters of $F_n$. As $l^2(F \backslash E) > 0$ and $F$ is not of full measure, there is a ball $B = B(x,r) \subset \Omega$ such that $l^2(B) > l^2(B \cap F) > 0$ and $\overline{B} \cap \overline{E} = \emptyset$. We notice that also for sufficiently large $n$ we have 

$$l^2(B) - c \geq l^2(B \cap F_n) \geq c > 0.$$
As all the sets $F_n$ have smooth boundary, we necessarily have $\partial F_n \cap B \neq \emptyset$. Let $y \in \overline{B}$ be the point closest to $\partial E$ and $l$ the line parallel to $\overline{p_1 p_2}$ passing through $y$. Finally, let $z_n \in \partial F_n \cap B$. The situation is presented on Figure \ref{fig:strictlyconvex}. 

\begin{figure}[h]
    \includegraphics[scale=0.27]{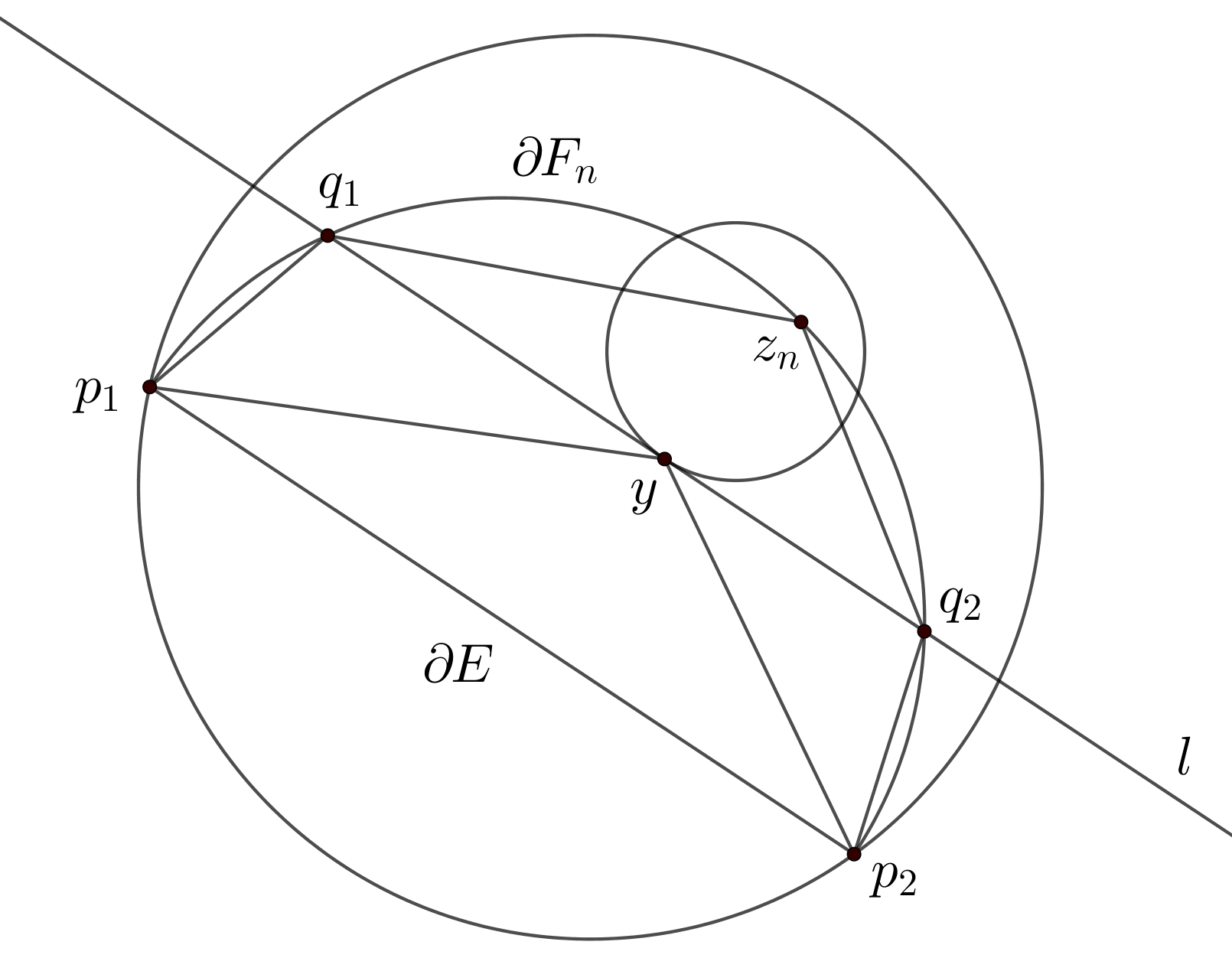}
    \caption{Construction of the lower bound for $P_\phi(F_n, \Omega)$}
    \label{fig:strictlyconvex}
\end{figure}
Let us denote by $\nu_1$ the vector normal to $\overline{p_1 z_n}$ and by $\nu_2$ the vector normal to $\overline{z_n p_2}$. By repeating the calculation from the proof of Proposition \ref{prop:linesegment} for the section of $\Gamma_n$ from $p_1$ to $z$ and then for the section of $\Gamma_n$ from $z_n$ to $p_2$, we see that we may estimate $P_\phi(F_n, \Omega)$ from below in the following way:

$$ P_\phi(F_n, \Omega) = \int_{\partial F_n} \phi(\nu(x)) \, d \mathcal{H}^{1} =\int_{\Gamma_n^{(p_1, z_n)}} \phi(\nu(x)) \, d \mathcal{H}^{1} + \int_{\Gamma_n^{(z_n, p_2)}} \phi(\nu(x)) \, d \mathcal{H}^{1} \geq $$
$$ \geq \phi(\nu_1) \mathcal{H}^1(\overline{p_1 z_n}) + \phi(\nu_2) \mathcal{H}^1(\overline{z_n p_2}),$$
Now, let $q_1 = \overline{p_1 z_n} \cap l$ and $q_2 = \overline{z_n p_2} \cap l$. Denote by $\nu_3$ be the vector normal to $\overline{p_1 y}$ and by $\nu_4$ the vector normal to $\overline{y p_2}$. By convexity of $\phi$ we have

$$ P_\phi(F_n, \Omega) \geq \phi(\nu_1) \mathcal{H}^1(\overline{p_1 z_n}) + \phi(\nu_2) \mathcal{H}^1(\overline{z_n p_2}) = \phi(\nu_1) \mathcal{H}^1(\overline{p_1 q_1}) + \phi(\nu_1) \mathcal{H}^1(\overline{q_1 z_n}) + $$
$$ + \phi(\nu_2) \mathcal{H}^1(\overline{z_n q_2}) + \phi(\nu_2) \mathcal{H}^1(\overline{q_2 p_2}) \geq \phi(\nu_1) \mathcal{H}^1(\overline{p_1 q_1}) + \phi(\nu_0) \mathcal{H}^1(\overline{q_1 y}) + \phi(\nu_0) \mathcal{H}^1(\overline{y q_2}) + $$
$$+ \phi(\nu_2) \mathcal{H}^1(\overline{q_2 p_2}) \geq \phi(\nu_3)\mathcal{H}^1(\overline{p_1 y}) + \phi(\nu_4)\mathcal{H}^1(\overline{y p_2}) = P.$$
By strict convexity of $\phi$ we have $P > \phi(\nu_0) \mathcal{H}^1(\overline{p_1 p_2}) = P_\phi(E, \Omega)$. As $F_n$ approximates $F$ is the strict topology with respect to $\phi$, we have

$$ P_\phi(F, \Omega) = \lim_{n \rightarrow \infty} P_\phi(F_n, \Omega) \geq P > P_\phi(E, \Omega).$$
But $E$ and $F$ have the same trace, by Proposition \ref{prop:linesegment} $E$ is a $\phi-$minimal set and we assumed $F$ to be a $\phi-$minimal set. Contradiction. Thus $E$ is the only $\phi$-minimal set with trace $T\chi_E$. \qed
\end{dd}

\begin{stw}\label{stw:regularityofminimalsets}
Suppose that $B_\phi(0,1)$ is strictly convex. Let $E$ be a $\phi-$minimal set. Then $\partial E = \bigcup_{i = 1}^\infty L_i$, where $L_i$ is a family of line segments, pairwise disjoint in $\overline{\Omega}$.
\end{stw}

\begin{dd}
Suppose that $x, y$ are two points in a path-connected component $S$ of $\partial^* E$. Then, by Proposition \ref{stw:strict}, the line segment $\overline{xy}$ lies in $S$. As the triangle inequality is strict, $S$ is a line segment; if there is another point $z \in S$, then again by Proposition \ref{stw:strict} $z$ is collinear with both $x$ and $y$. Thus every path-connected component of $\partial^* E$ is a line segment. 

As (with our choice of representative) $\partial^* E$ is dense in $\partial E$, we have $\partial^* E = \partial E$; take $x \in \partial E \backslash \partial^* E$. Then in some neigbhourhood of $x$ there are infinitely many connected components of $\partial^* E$; in particular, the $\phi-$perimeter of $E$ is infinite, contradiction. Thus $\partial E = \bigcup_{i = 1}^\infty L_i$, where $L_i$ are (pairwise disjoint) line segments. \qed
\end{dd}

But this leads directly to uniqueness of solutions for continuous boundary data; as we know that the only $\phi$-minimal surfaces are line segments, the proof of \cite[Theorem 4.1]{SWZ} holds with minimal changes and we have:

\begin{wn}
Suppose that $\Omega$ is strictly convex and that $B_\phi(0,1)$ is strictly convex. Then the solution to problem (\ref{problem}) with continuous boundary data is unique. \qed
\end{wn}

Now, we turn to the issue of existence of solutions to problem (\ref{problem}). To this end, recall the barrier condition (Definition \ref{def:barrier}). By \cite[Theorem 1.1]{JMN} it is a sufficient condition for existence of solutions for continuous boundary data. However, in dimension two, Proposition \ref{stw:strict} implies that the class of sets that satisfy the barrier condition are precisely the open bounded strictly convex sets (regardless of the choice of $\phi$). Thus, if $\Omega \subset \mathbb{R}^2$ is an open bounded strictly convex set and $B_\phi(0,1)$ is strictly convex, then there exists a solution to the least gradient problem. We summarize the above discussion in

\begin{tw}\label{tw:strict}
Let $\Omega \subset \mathbb{R}^2$ be an open bounded strictly convex set. Suppose that $\phi$ is a metric integrand, $\phi(x, Du) = \phi(Du)$ and $B_\phi(0,1)$ is strictly convex. Let $f \in C(\partial\Omega)$. Then then there exists a unique solution to Problem (\ref{problem}). \qed
\end{tw}

\begin{uw}
While the existence proof is based on \cite[Theorem 1.1]{JMN}, let us note that the above result is substantially different from the uniqueness theorem \cite[Theorem 1.2]{JMN}. That result requires $\phi$ to satisfy two additional conditions: firstly, that $\phi$ has regularity somewhat stronger than $W^{2,\infty}$ outside the origin; what is more relevant in our case, the other assumption states that there exists $C > 0$ such that

$$ \sum_{i,j = 1}^n \phi_{\xi_i \xi_j} (x, \xi) p^i p^j \geq C |p - (p \cdotp \xi) \xi|^2$$
for every $p \in \mathbb{R}^n$ and every $\xi \in S^{n-1}$.
This is a convexity assumption stronger than strict convexity of $B_\phi(0,1)$. Thus the results proved in the this Section, and even more so in the next, are independent from the results in \cite{JMN}: the case when $B_\phi(0,1)$ is not strictly convex is not covered at all by \cite[Theorem 1.2]{JMN}, while in the case when $B_\phi(0,1)$ is strictly convex we are able to prove uniqueness of minimizers to Problem (\ref{problem}) regardless of the regularity or stronger modes of convexity of $\phi$. However, we note that in our setting we only allow $\phi$ to depend on the direction of the derivative and not on location.
\end{uw}

Here we also recall three useful results, proved in \cite[Proposition 3.5]{Gor1}, \cite[Lemma 2.8]{Gor2} and \cite[Lemma 3.8]{GRS} respectively in the isotropic case.

\begin{stw}
Suppose that $\Omega$ is convex, $B_\phi(0,1)$ is strictly convex and suppose that $u \in BV(\Omega)$ is a function of $\phi-$least gradient. Then for every $t \in \mathbb{R}$ we have $\partial \{ u > t \} = \bigcup_{i=1}^\infty L_{t,i}$, where $L_{t,i}$ are line segments with ends on $\partial\Omega$ and this union is locally finite in $\Omega$. Furthermore, $\overline{L_{t,i}}$ are pairwise disjoint in $\overline{\Omega}$.
\end{stw}

The proof in the isotropic case relies only on the regularity of $\phi-$minimal sets and the fact that the triangle inequality is strict. Both of these facts are true if $B_\phi(0,1)$ are strictly convex, so the Proposition remains true.

\begin{lem}\label{lem:wewnatrz}
Suppose that $B_\phi(0,1)$ is strictly convex and that $u \in BV(\Omega)$ is a function of $\phi-$least gradient. Let $E_t = \{ u \geq t \}$. Suppose that $x \in \Omega$ is a point of continuity of $u$, $u(x) = t$ and $x \notin \partial E_t$. Then there exists a ball $B(x,r) \subset E_t$. 
\end{lem}

The proof in the isotropic case relies only on the regularity of $\phi-$minimal sets and the relative isoperimetric inequality. Both of these facts hold for strictly convex $B_\phi(0,1)$, so the Lemma remains true.

\begin{lem}\label{lem:grs}
Suppose that $\Omega$ is strictly convex and $B_\phi(0,1)$ is strictly convex and suppose that $u$ is a minimizer of Problem (\ref{problem}). Let $f \in C(\partial\Omega)$. Then for every $t \in \mathbb{R}$ we have
$$ \partial \{ u \geq t \} \cap \partial\Omega \subset f^{-1}(t).$$
\end{lem}

The proof in the isotropic case relies only on the regularity of $\phi-$minimal sets and a blow-up argument, which is applied at regular points of $\partial\Omega$ and does not depend on $\phi$. Thus the Lemma remains true in the anisotropic case.

\subsection{Regularity of minimizers}

We briefly recall the regularity results from \cite{SWZ} concerning the isotropic case. The authors assume that $\partial \Omega$ is of class $C^2$ and that $\Omega$ is uniformly convex, i.e. the mean curvature of $\partial\Omega$ is positive. Then if the boundary data $f$ is of class $C^{0,\alpha}(\partial\Omega)$, where $\alpha \in (0,1]$, then the corresponding minimizer to the least gradient problem $u$ is in the class $C^{0, \alpha/2}(\Omega)$. A similar result is obtained if the mean curvature can vanish at isolated points and has polynomial growth. As the (two-dimensional) examples provided by the authors show, the above results are optimal.

Our goal is to extend these results to the anisotropic case. This issue has been recently explored in \cite{DS}; the authors use a different approach, going through the optimal transport theory and using the equivalence proved in \cite{GRS}, and are able to prove regularity estimates for $W^{1,p}$ boundary data (with $p \leq 2$). Here, we discuss the issue of regularity of minimizers when the boundary data are not weakly differentiable, expressed only in terms of the modulus of continuity of the boundary data. Moreover, as we use a different approach, we may replace the regularity assumptions on the $\partial\Omega$ by weaker ones. However, \cite[Example 5.8]{SWZ} shows that we cannot get rid of some form of uniform convexity altogether. In this paper, we will allow $\Omega$ to have only Lipschitz boundary and use the following definition of uniform convexity (which agrees with the classical definition for $C^2$ sets, see Proposition \ref{stw:uniformly}):

\begin{dfn}\label{dfn:uniformlyconvex}
We say that the set an open bounded convex set $\Omega$ is uniformly convex, if the following condition is satisfied: let $P = \{ y \geq a x^2 \}$, where $a > 0$. Let $x_0 \in \partial\Omega$ and let $l$ be a supporting line at $x_0$. Then there exists $P'$, an isometric image of $P$, tangent to $l$ at $x_0$ such that $\overline{\Omega} \subset P'$ and $\partial P' \cap \overline{\Omega} = \{ x_0 \}$.
\end{dfn}

Similarly, we will say that an open bounded set $\Omega$ is $\beta-$uniformly convex, if for some $\beta > 0$ the Definition above is satisfied with $\widetilde{P} = \{ y \geq a x^{\beta + 2} \}$ in place of $P$.

\begin{stw}\label{stw:uniformly}
If $\partial \Omega \in C^2$ and its curvature is positive, then $\Omega$ is uniformly convex in the sense of Definition \ref{dfn:uniformlyconvex}.
\end{stw}

\begin{dd}
As $\partial\Omega$ is compact, the mean curvature has a positive lower bound $c$. Take any $x \in \partial\Omega$ and let $l$ be a line tangent to $\partial\Omega$ at $x$. We choose the coordinate system so that $x = (0,0)$ and $l = \{ y = 0 \}$. As $\partial\Omega$ is strictly convex, it is a union of two graphs of convex functions. Let $g$ be one of these functions and $g(0) = 0$. Then, by the formula for the curvature of a graph, we have
$$ g'' = (1 + (g')^2)^{3/2} k \geq k \geq c,$$
so we have $g(0) = 0$, $g'(0) = 0$ and $g'' \geq c$. Thus $g(x) \geq \frac{c}{2} x^2$. As for the second function, its graph lies above the graph of $g$, so also above the parabola $y = \frac{c}{2} x^2$. As the coefficient does not depend on $x \in \partial\Omega$, $\Omega$ is uniformly convex in the sense of Definition \ref{dfn:uniformlyconvex}. \qed
\end{dd}

With essentially the same proof we obtain

\begin{wn}
If $\partial\Omega \in C^2$ and for every $x_0 \in \partial\Omega$ its curvature satisfies a bound $k(x) \geq a|x - x_0|^\beta$ in some neighbourhood of $x_0$, then $\Omega$ is $\beta-$uniformly convex. \qed
\end{wn}

\begin{uw}
The condition in Definition \ref{dfn:uniformlyconvex} seems as if it was hard to check for any given set $\Omega$. However, it is sufficient to check it for every $x_0$ for at most two supporting lines: without loss of generality assume that $x_0 = (0,0)$. The set of supporting lines $y = ax$, parametrized by the coefficient $a$, is closed and convex, so it is an interval $[a_1, a_2]$. Take the supporting parabolas $P_1$ and $P_2$ corresponding to lines $y= a_1 x$ and $y = a_2 x$. Take the parabola $P$ corresponding to a line $y = ax$, where $a \in (a_1, a_2)$. Then $P_1 \cap P_2 \subset P$ and in particular $\Omega \subset P_1 \cap P_2 \subset P$; this can be easily seen in the polar form of the equation for the parabola. Thus, if the condition from the definition of uniform convexity is satisfied for two extreme supporting lines at $x_0$, it is satisfied for all supporting lines at $x_0$.

In particular, an important class of sets $\Omega$ uniformly convex in the sense of Definition \ref{dfn:uniformlyconvex} are strictly convex sets such that $\partial\Omega$ is $C^2$ except for finitely many corners and the curvature of $\Omega$ is bounded from below (on the set where $\partial\Omega$ is $C^2$).
\end{uw}

We turn our attention to the regularity of minimizers to Problem (\ref{problem}). Firstly, we use a variant of a result from \cite{JMN} to prove that any minimizer is continuous up to the boundary of $\Omega$. We recall that \cite[Theorem 1.3]{JMN}, which asserts the continuity of solutions up to the boundary, follows from the following comparison principle:

\begin{stw}
(\cite[Theorem 4.6]{JMN}) Let $\Omega \subset \mathbb{R}^2$ be an bounded convex set. If $\Omega$ satisfies the barrier condition, $E_1, E_2 \subset \mathbb{R}^2$ are $\phi-$area minimizing in $\Omega$ and

$$ E_1 \backslash \Omega \subset \subset E_2 \backslash \Omega,$$
then $E_1 \subset \subset E_2$.
\end{stw}
While this result is originally stated without restrictions as to the dimension, but for a metric integrand satisfying some additional regularity properties, as we have a special form of $\phi$, i.e. it depends only on the direction of derivative of $u$, the above Proposition follows directly from Proposition \ref{stw:regularityofminimalsets}. Thus we obtain

\begin{wn}
Suppose that $\Omega \subset \mathbb{R}^2$ is an open bounded convex set. Let $\phi$ be a metric integrand such that $B_\phi(0,1)$ is strictly convex and $u$ is a solution of Problem (\ref{problem}) with boundary data $f$. Then $u \in C(\overline{\Omega})$.
\end{wn}

We recall that an increasing function $\omega: [0,\infty] \rightarrow [0,\infty]$ such that $\omega(0^+) = 0$ is a modulus of continuity of a continuous function $f$, if $|f(x) - f(y)| \leq \omega(|x-y|)$. The next Proposition is our main regularity result; we present it for uniformly convex sets for the sake of clarity and then show how to pass to $\beta-$uniformly convex sets.

\begin{stw}\label{prop:modulusofcontinuity}
Suppose that $\Omega \subset \mathbb{R}^2$ is uniformly convex and $B_\phi(0,1)$ is strictly convex. Let $f \in C(\partial\Omega)$ and take $\omega$ to be its modulus of continuity. Let $u$ be the solution of Problem (\ref{problem}) with boundary data $f$. Then $u \in C(\overline{\Omega})$ and it is continuous with modulus of continuity

$$ \overline{\omega}(|x - y|) = \omega(c(\Omega) |x - y|^{1/2}).$$
\end{stw}

\begin{dd}
The proof will follow in three steps. Firstly, we prove the statement in a special geometric situation and then gradually reduce the general case to the special case.

Step 1. Let $p, q \in \Omega$. Suppose that $p \in \partial E_t$ and $q \in \partial E_s$. Let $l_p$ be a line passing through $p$ such that the connected component of $\partial E_t$ containing $p$ lies inside $l_p$ (and similarly we define $l_q$). Suppose that $l_p$ and $l_q$ are parallel and the line $pq$ is perpendicular to $l_p$ (and $l_q$). Let $x_0 \in \partial \Omega$ be a point such that there is a supporting line $l$ at $x_0$ parallel to $l_p$ (and $l_q$); there are two such points, without loss of generality $l_p$ is closer to $x_0$ then $l_q$. 

We change coordinates so that $x_0 = (0,0)$, $l = \{ y = 0 \}$ and $x_0$ is the lowest point of $\overline{\Omega}$. Take a supporting parabola $P = \{ y = a x^2 \}$ at $x_0$ as in Definition \ref{dfn:uniformlyconvex}. Let $p' \in \partial \Omega \cap l_p$ and $p'' \in \partial P \cap l_p$. Similarly we define $q'$ and $q''$ and we require that $p'$ and $q'$ (and $p''$ and $q''$) lie on the same side with respect to the vertical line $pq$. The situation is presented on Figure \ref{fig:regularity}.

\begin{figure}[h]
    \includegraphics[scale=0.27]{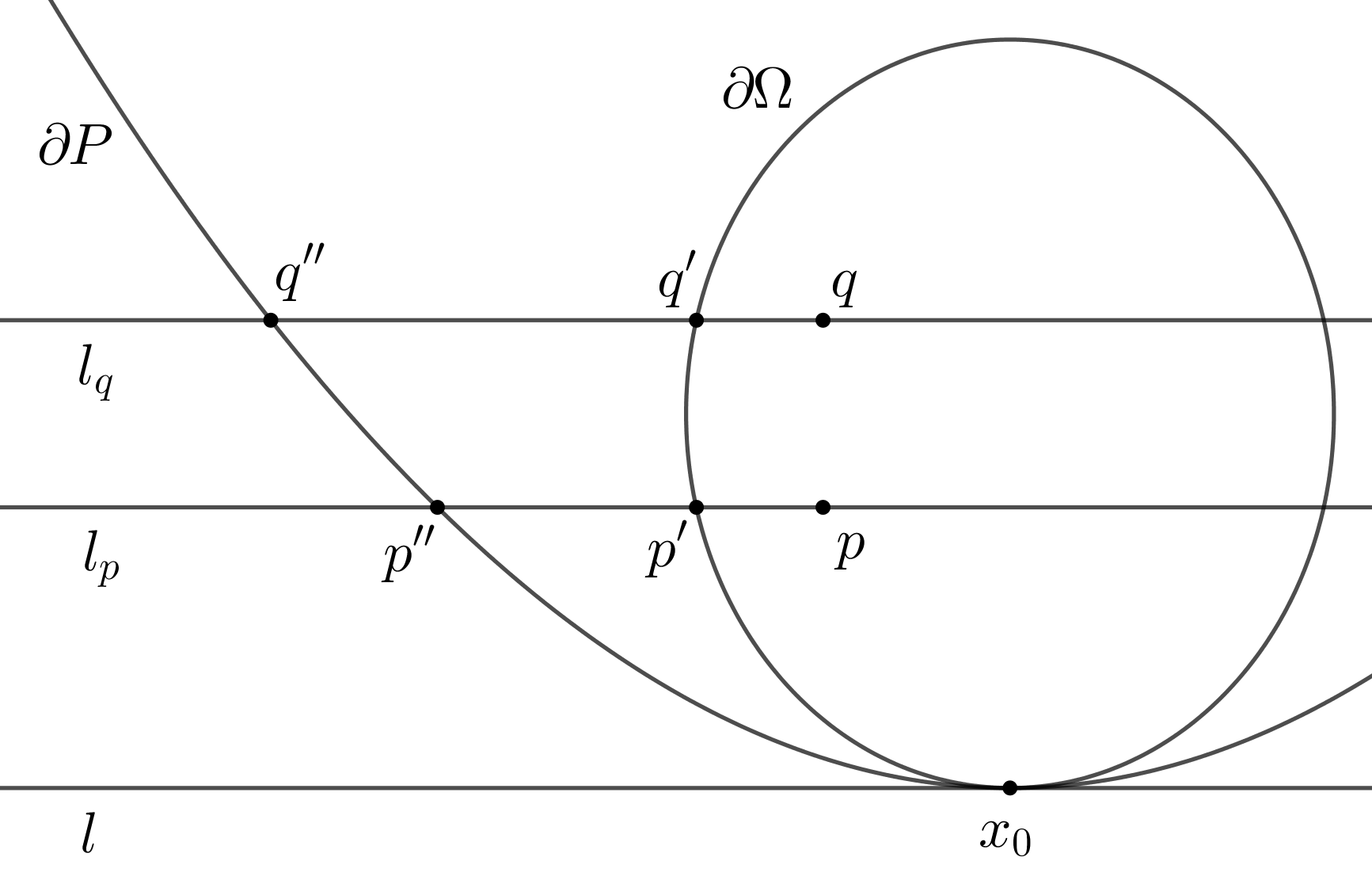}
    \caption{The construction of the competitor $F$}
    \label{fig:regularity}
\end{figure}

Since $u \in C(\overline{\Omega})$, we have

$$ |u(q) - u(p)| = |u(q') - u(p')| = |f(q') - f(p')| \leq \omega(|q' - p'|). $$
Now, we read off from the geometrical situation that $|q' - p'| \leq |q'' - p''|$: take a supporting line at $p'$. Take a parabola $P'$, isometric to $P$ such as in Definition \ref{dfn:uniformlyconvex}. Then, by Definition \ref{dfn:uniformlyconvex}, $|q' - p'|$ is smaller than $|\widetilde{q'} - p'|$, where $\widetilde{q'}$ lies on the intersection of $P'$ and $l_q$. As $P'$ is an isometric image of $P$, the curvature of $P'$ at $p'$ is bounded from below by the curvature of the parabola $P + p'$ at $p'$, then as in the proof of Proposition 5.2 while $|q_y'' - p_y''| = |q_y' - p_y'|$, we have $|q_x'' - p_x''| > |q_x' - p_x'|$. 

Furthermore, as $\partial P$ is a parabola, we see that $|q'' - p''|^2 \leq C(\Omega) |q - p|$. We calculate

$$ |q'' - p''|^2 = (q_y'' - p_y'')^2 + (q_x'' - p_x'')^2 = d(l_p, l_q)^2 + (\sqrt{\frac{d(l_q, l)}{a}} - \sqrt{\frac{d(l_p, l)}{a}})^2 = $$
$$ = d(l_p, l_q)^2 + \frac{d(l_q,l)}{a} + \frac{d(l_p, l)}{a} - \frac{2}{a} \sqrt{d(l_q, l) d(l_p, l)} = d(l_p, l_q)^2 + \frac{d(l_q,l_p)}{a} + $$
$$ + \frac{2}{a} (d(l_p, l) - \sqrt{d(l_q, l) d(l_p, l)}) \leq d(l_p, l_q) (d(l_p, l_q) + \frac{1}{a}) + 0 \leq (\text{diam }\Omega + \frac{1}{a}) d(l_p, l_q) = c(\Omega) |q - p|.$$
Returning to the level of moduli of continuity, we obtain that

$$ |u(q) - u(p)| \leq \omega(|q' - p'|) \leq \omega(|q'' - p''|) \leq \omega(c(\Omega) |q - p|^{1/2}) = \overline{\omega}(|q - p|).$$

Step 2. Let $p, q \in \Omega$ and suppose that $p \in \partial E_t$ and $q \in \partial E_s$. We prove that it is sufficient to assume that $p$ and $q$ are as assumed in Step 1. We use the notation as in Step 1, except now $l_p$ is not parallel to $l_q$ (so $l$ is parallel only to $l_p$ and not to $l_q$). Since $u \in C(\overline{\Omega})$, $l_p$ and $l_q$ intersect somewhere outside of $\overline{\Omega}$. Let $p'$ and $q'$ be on this side of the line $pq$ so that $q_y'$ is smaller. Now, we draw a line $l_q'$ parallel to $l_p$ passing through $q'$. By Step 2 we have

$$ |u(q) - u(p)| = |f(q') - f(p')| \leq \overline{\omega}(d(l_p, l_q')) \leq \overline{\omega}(|q - p|).$$

Step 3. Finally, we prove that it is sufficient to assume that $p$ and $q$ are as assumed in Step 2. Firstly, we notice that if $p$ lies on $\partial \{ u > t \}$ instead of $\partial E_t$, then the proofs in Steps 1 and 2 remain unchanged. Let $p, q \in \Omega$. If $u(p) = u(q)$, then there is nothing to prove; without loss of generality $u(p) > u(q)$. Let $u(p) = t$ and suppose that $p \notin \partial E_t$ and $p \notin \partial \{ u > t \}$. By Lemma \ref{lem:wewnatrz} $p$ lies on a level set $U$ of $u$ of a positive measure. Then we may replace $p$ with $\widetilde{p}$, which lies on $\partial U$, so on $\partial E_t$ or $\partial \{ u > t \}$ and is closer to $q$ than $p$. A similar analysis applies to $q$. \qed
\end{dd}

Let us stress the fact that in the Proposition above the constant $c(\Omega)$ depends only on $\Omega$ and not on the metric integrand $\phi$ (and is explicitly given in the calculation in Step 1). If the boundary is $C^2$, then it depends only on the diameter and the lower bound $c$ on curvature of $\Omega$: by Step 1 and Proposition \ref{stw:uniformly} we have $c(\Omega) =$ diam $\Omega + \frac{1}{a} =$ diam $\Omega + \frac{2}{c}$.

\begin{wn}
Under the assumptions of Proposition \ref{prop:modulusofcontinuity}, suppose that $f \in C^{0, \alpha}(\partial\Omega)$. Then $u \in C^{0, \alpha/2}(\overline{\Omega})$. \qed
\end{wn}

Now, we will formulate a similar result for $\beta-$uniformly convex sets.

\begin{stw}\label{prop:modulusofbetacontinuity}
Suppose that $\Omega$ is $\beta-$uniformly convex and $B_\phi(0,1)$ is strictly convex. Let $f \in C(\partial\Omega)$ and take $\omega$ to be its modulus of continuity. Let $u$ be the solution of Problem (\ref{problem}) with boundary data $f$. Then $u \in C(\overline{\Omega})$ and it is continuous with modulus of continuity

$$ \overline{\omega}(|x - y|) = \omega(c(\Omega) |x - y|^{1/(\beta + 2)}).$$
\end{stw}

\begin{dd}
The proof remains unchanged except for the final calculation of distances in Step 1, where we obtain
$$ |q'' - p''|^2 = (q_y'' - p_y'')^2 + (q_x'' - p_x'')^2 = d(l_p, l_q)^2 + (({\frac{d(l_q, l)}{a}})^{\frac{1}{\beta + 2}} - ({\frac{d(l_p, l)}{a}})^{\frac{1}{\beta + 2}})^2 = $$
$$ = d(l_p, l_q)^2 + (\frac{d(l_q,l)}{a})^{\frac{2}{\beta + 2}} + (\frac{d(l_p, l)}{a})^{\frac{2}{\beta + 2}} - 2 (\frac{d(l_q, l) d(l_p, l)}{a^2})^{\frac{1}{\beta + 2}} \leq d(l_p, l_q)^2 + (\frac{d(l_q,l_p)}{a})^{\frac{2}{\beta + 2}} $$
$$ + (\frac{d(l_p, l)}{a})^{\frac{2}{\beta + 2}} + (\frac{d(l_p, l)}{a})^{\frac{2}{\beta + 2}} - 2 (\frac{d(l_q, l) d(l_p, l)}{a^2})^{\frac{1}{\beta + 2}} \leq d(l_p, l_q)^{\frac{2}{\beta + 2}} (d(l_p, l_q)^{2 - \frac{2}{\beta + 2}} + (\frac{1}{a})^{\frac{2}{\beta+2}}) + 0 \leq $$
$$\leq ((\text{diam }\Omega)^{2 - \frac{2}{\beta + 2}} + (\frac{1}{a})^{\frac{2}{\beta + 2}}) d(l_p, l_q)^{\frac{2}{\beta + 2}} = c(\Omega) |q - p|^{\frac{2}{\beta + 2}}.$$

\end{dd}

\begin{wn}
Under the assumptions of Proposition \ref{prop:modulusofbetacontinuity}, suppose that $f \in C^{0, \alpha}(\partial\Omega)$. Then $u \in C^{0, \alpha/(\beta + 2)}(\overline{\Omega})$. \qed
\end{wn}

We conclude with a comparison of the above results with the results in \cite[Section 5]{SWZ}. We see that in dimension two we obtained the same regularity estimates as in the anisotropic case. As the (counter)examples in \cite{SWZ} show, these results are optimal. Moreover, this line of reasoning enables us to prove them with weaker assumptions concerning the regularity of $\partial\Omega$.

\section{Non-strictly convex unit ball $B_\phi(0,1)$}

When $\phi$ is a metric integrand such that $B_\phi(0,1)$ is not strictly convex, then no open set $\Omega$ with $C^1$ boundary satisfies the barrier condition (the proof of this fact is presented at the end of this Section). Thus, existence of minimizers is not guaranteed and we have to prove it by using another means. Furthermore, we will see in Lemma \ref{lem:nonuniqueminimalsets} that line segments are not the only connected $\phi-$minimal surfaces; thus we may not use the reasoning from \cite{SWZ} to conclude uniqueness of minimizers. Throughout this Section, $B_\phi(0,1)$ is convex but not strictly convex and $I$ always denotes a line segment in $\partial B_\phi(0,1)$.

This Section is organized as follows. Firstly, we see that if $\Omega$ is uniformly convex (or $\beta-$uniformly convex), then there exists a solution to Problem (\ref{problem}) and it satisfies the regularity estimates as proved for strictly convex $B_\phi(0,1)$ in Proposition \ref{prop:modulusofcontinuity}. Furthermore, we prove that if $\Omega$ is only strictly convex, there still exists a minimizer of Problem (\ref{problem}). Secondly, we will show that line segments are not the only $\phi-$minimal surfaces and infer that minimizers to Problem (\ref{problem}) may fail to be unique even for smooth boundary data (however, there still may exist boundary data for which minimizers are unique; see \cite[Example 5.15]{Gor1}). Thirdly, we show that there exist such boundary data $f \in C^\infty(\partial\Omega)$ such that some minimizers of Problem (\ref{problem}) have regularity no better than $BV(\Omega) \cap L^\infty(\Omega)$. Finally, we will see why the barrier condition is not satisfied for sets with $C^1$ boundary.

\begin{tw}\label{tw:existence}
Let $\Omega \subset \mathbb{R}^2$ be an open bounded uniformly convex set. Suppose that $\phi$ is a metric integrand depending only on the second variable. Let $f \in C(\partial\Omega)$. Then there exists a solution $u \in C(\overline{\Omega})$ to Problem (\ref{problem}). Additionally, if $\omega$ is the modulus of continuity of $f$, then $\overline{\omega}$ defined as in Proposition \ref{prop:modulusofcontinuity} is the modulus of continuity of $u$.
\end{tw}

\begin{dd}
Let $\omega$ be a modulus of continuity of $f$ (as $\partial \Omega$ is compact, $f$ is uniformly continuous, so it admits a modulus of continuity). Take $\phi$ to be any metric integrand. Then $\phi_n = \phi + \frac{1}{n} l^2$ is such that $B_{\phi_n}(0,1)$ is strictly convex. By Theorem \ref{tw:strict} there exists a solution $u_n \in C(\overline{\Omega})$ to Problem (\ref{problem}) with boundary data $f$ with respect to the anisotropic norm $\phi_n$.

By Proposition \ref{prop:modulusofcontinuity} the solution $u_n$ is continuous on $\overline{\Omega}$ with modulus of continuity $\overline{\omega}$, which depends only on the geometry of $\Omega$ and not on the metric integrand $\phi_n$; thus the sequence $u_n$ has the same modulus of continuity, so it is equicontinuous. Also, the sequence $u_n$ is uniformly bounded from below and above by the maximum and minimum of $f$. By Arzela-Ascoli theorem the sequence $u_n$ admits a subsequence which converges uniformly on the compact set $\overline{\Omega}$.

We obtain that $u_{n_k} \rightarrow u$ uniformly on $\overline{\Omega}$. In particular, this convergence is in $L^1(\Omega)$, hence $u \in BV(\Omega)$, so also $u \in BV_\phi(\Omega)$. As the convergence is uniform, we have $Tu = u|_{\partial\Omega} = f$. By Theorem \ref{tw:gammaconvergence} $u$ is a function of $\phi-$least gradient (note that we may not use Miranda's theorem, since we change the anisotropic norm). Finally, we observe that as each $u_n$ admitted the same modulus of continuity $\overline{\omega}$, $u$ is uniformly continuous with the same modulus of continuity. \qed
\end{dd}

An analogous proof shows that the above result also holds if $\Omega$ is only $\beta-$uniformly convex. In the next result, we show that if $\Omega$ is only strictly convex, we still obtain existence of minimizers for continuous boundary data; however, in this case we do not have any regularity estimates.

\begin{tw}\label{tw:existence2}
Let $\Omega \subset \mathbb{R}^2$ be an open bounded strictly convex set. Suppose that $\phi$ is a metric integrand depending only on the second variable. Let $f \in C(\partial\Omega)$. Then there exists at least one solution to Problem (\ref{problem}).
\end{tw}

\begin{dd}
1. Due to Theorem \ref{tw:strict} we only have to prove the result if $B_\phi(0,1)$ is not strictly convex. Then $\phi_n = \phi + \frac{1}{n} l^2$ is such that $B_{\phi_n}(0,1)$ is strictly convex. By Theorem \ref{tw:strict} there exists a solution $u_n \in C(\overline{\Omega})$ to Problem (\ref{problem}) with boundary data $f$ with respect to the anisotropic norm $\phi_n$. Furthermore, the family $u_n$ is uniformly bounded in $BV(\Omega)$, as

$$ \int_\Omega |u_n|dx + \int_\Omega |Du_n| \leq \int_\Omega \sup_{\partial\Omega} |f| dx + C \int_\Omega |Du_n|_\phi + C \int_\Omega \frac{1}{n} |Du_n| \leq |\Omega| \sup_{\partial\Omega} |f| + C F_{\phi_n}(u_n) \leq$$
$$ \leq |\Omega| \sup_{\partial\Omega} |f| + C F_{\phi_n}(v \equiv 0) \leq |\Omega| \sup_{\partial\Omega} |f| + C \int_{\partial\Omega} (\sup_{\partial B(0,1)} \phi + \frac{1}{n}) |f| d\mathcal{H}^{1} \leq M.$$
In particular, $u_n$ admits a subsequence (still denoted by $u_n$) convergent in $L^1(\Omega)$. By Theorem \ref{tw:gammaconvergence} $u$ is a minimizer of the functional $F_\phi$; if we prove additionally that $Tu = f$, then $u$ is a solution to Problem (\ref{problem}).

2. We recall that if the trace of $u$ equals $f$, then the set $T$ of such $x \in \partial\Omega$ that
$$ \dashint_{B(x,r) \Omega} |u(y) - f(x)| dy \rightarrow 0$$
when $r \rightarrow 0$ is of $\mathcal{H}^1-$full measure (see \cite[Theorem 5.3.2]{EG}). Fix $x \in T$ and arbitrary $\varepsilon > 0$. As $f \in C(\partial\Omega)$, there exists a neighbourhood of $x$ in $\partial\Omega$ such that 
$$ f(x) - \varepsilon \leq f(y) \leq f(x) + \varepsilon \qquad \text{ in } B(x, \delta_1) \cap \partial\Omega.$$
3. As $\Omega$ is strictly convex, for sufficiently small $\delta_1$ the set $B(x, \delta_1) \cap \partial\Omega$ consists of two points $p_1, p_2$ and the line segment $\overline{p_1 p_2}$ lies inside $\Omega$. Denote by $\Delta$ the open set bounded by an arc of $\partial\Omega$ containing $x$ and the line segment $\overline{p_1 p_2}$. Let us take a ball $B(x, \delta_2)$ such that $B(x, \delta_2) \cap \Omega \subset \Delta$. Then for every $n$ we have
$$ f(x) - \varepsilon \leq u_n(y) \leq f(x) + \varepsilon \qquad \text{ in } B(x, \delta_2) \cap \Omega;$$
suppose otherwise, i.e. that for some $y \in B(x, \delta_2) \cap \Omega$ we have $y \in \partial \{u_n \geq t \}$, where $t > f(x) + \varepsilon$. Take the connected component $S$ of $\partial \{ u_n \geq t \}$ containing $y$. By Proposition \ref{lem:grs} we have $S \cap \partial\Omega = \{ q_1, q_2 \} \subset f^{-1}(t)$. As $y \in B(x, \delta_2) \cap \Omega \subset \Delta$, at least one of points $q_1, q_2$ lies on $\partial\Omega \cap \partial\Delta$; on the other hand, by Step 2 we have $u_n(y) \leq f(x) + \varepsilon < t$ on $\partial \Omega \cap \partial \Delta$, contradiction. The case when $t < f(x) - \varepsilon$ is handed similarly.

4. As $u_n \rightarrow u$ in $L^1(\Omega)$, on some subsequence (still denoted $u_n$) we have convergence almost everywhere; hence
$$ f(x) - \varepsilon \leq u(y) \leq f(x) + \varepsilon \qquad \text{ for a.e. } y \in B(x, \delta_2) \cap \Omega.$$

5. Now, suppose that $Tu(x) = a > f(x)$. As $\varepsilon$ was arbitrary, we choose it to be small enough to satisfy $a > f(x) + \varepsilon$. Then for $r < \delta_2(\varepsilon)$

$$ \dashint_{B(x,r) \cap \Omega} |u(y) - a| dy = \dashint_{B(x,r) \cap \Omega} |u(y) - (f(x) + \varepsilon) + (f(x) + \varepsilon - a)| dy = $$
$$ = \dashint_{B(x,r) \cap \Omega} |u(y) - (f(x) + \varepsilon)| + \dashint_{B(x,r) \cap \Omega} |(f(x) + \varepsilon - a)| dy \geq 0 + |(f(x) + \varepsilon - a)|,$$
but as $x \in T$, the mean integral should vanish in the limit $r \rightarrow 0$, contradiction. A similar argument covers the case when $Tu(x) < f(x)$. Thus $Tu(x) = f(x)$ almost everywhere with respect to $\mathcal{H}^1$, so $Tu = f$. \qed

\end{dd}

We turn our attention to the issue of uniqueness of solutions. The key idea here is that we may perturb the level sets of a solution as long, as {\it a positive multiple of a normal vector $\alpha \nu \in I \subset \partial B_\phi(0,1)$}. 

\begin{uw}\label{uw:neighbourhood}
We notice that if $\alpha \nu \in \text{int } I$, then there exists a neighbourhood $N \subset S^1$ of $\nu_0$ such that for each $\nu \in N$ a positive multiple of $\nu$, namely $(\nu_0 \cdotp \nu)^{-1} \alpha \nu$, lies in $I$. 
\end{uw}

We use this observation to construct $\phi-$minimal surfaces other than a line segment. The next result is a construction of a $\phi-$minimal set with boundary which is not a line segment; later, we are going to extensively use this construction to prove existence of solutions with regularity no better than $BV(\Omega) \cap L^\infty(\Omega)$.

\begin{lem}\label{lem:nonuniqueminimalsets}
Suppose that there is a line segment $I \subset \partial B_\phi(0,1)$. Let $p_1, p_2 \in \partial\Omega$, take $\nu_0$ to be a vector normal to $\overline{p_1 p_2}$ and suppose that $\alpha \nu_0 \in \text{int } I$. Let $E \subset \Omega$ be an open set such that its boundary is the line segment $p_1 p_2$. Let $F \subset \Omega$ be an open set such that its boundary is a (finite) polygonal chain $p_1 q_1 ... q_n p_2$ such that the normal vector to each of the line segments in this polygonal chain lies in $N$.
Then

$$ P_\phi(F, \Omega) = P_\phi(E, \Omega).$$
In particular, in the notation of Proposition \ref{prop:linesegment}, the set $E$ is not the only $\phi-$minimal set with the same boundary data.
\end{lem}

\begin{dd}
Rename the points so that $q_0 = p_1$ and $q_{n+1} = p_2$. Take the polygonal chain $q_0 q_1 ... q_n q_{n+1}$ as in the assumption of the Lemma. It is enough to show that the set $F$ bounded by this polygonal chain has the same anisotropic perimeter as the set $F'$ bounded by the polygonal chain $q_0 q_2 ... q_n p_{n+1}$; then, we use this result to expand our line segment into a polygonal chain using finitely many steps without changing the perimeter.

\begin{figure}[h]
    \includegraphics[scale=0.27]{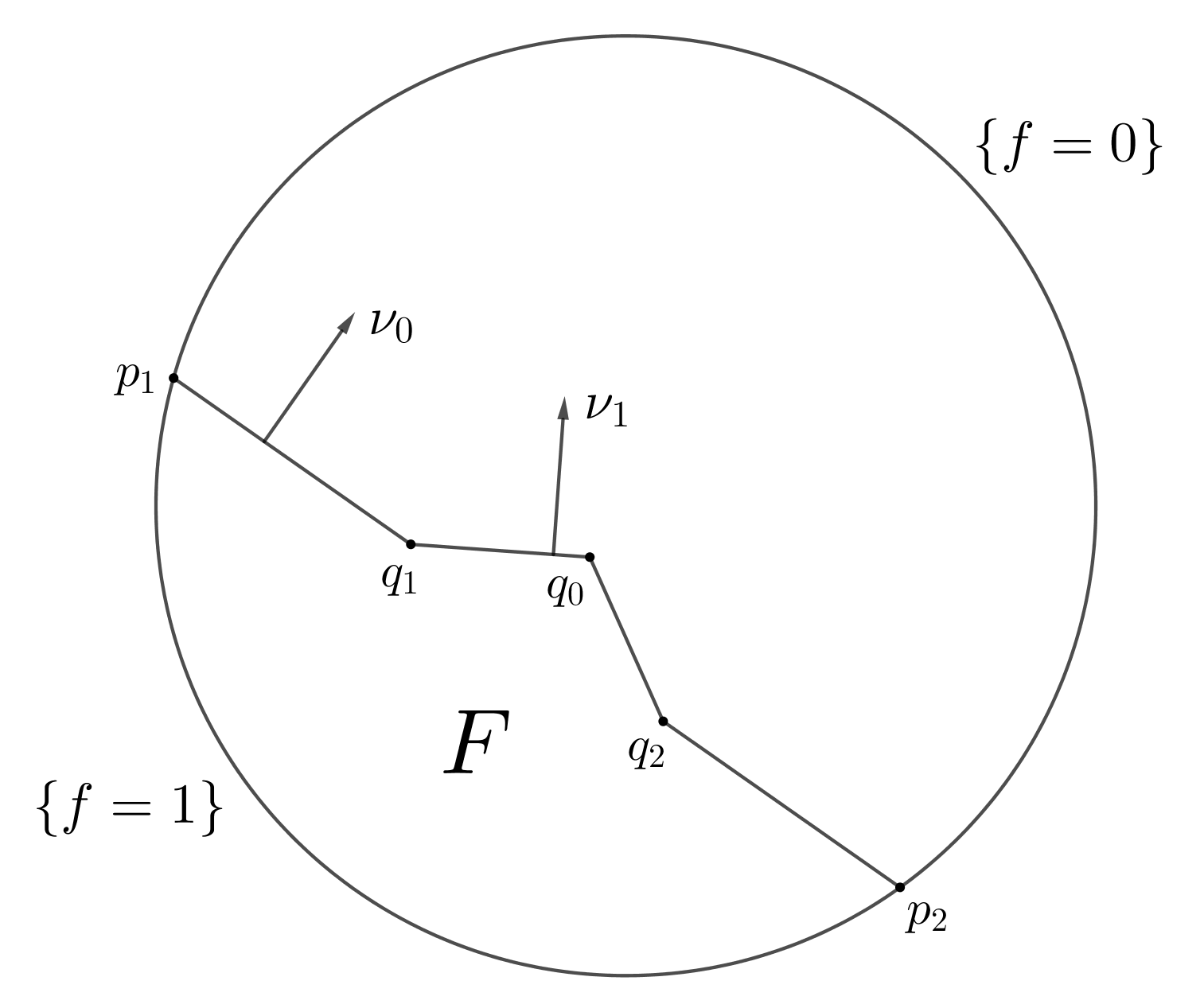}
    \caption{The construction of the competitor $F$}
    \label{fig:odcinki}
\end{figure}

We calculate the perimeter of $F$. As $\partial F$ is piecewise $C^1$ (and the measure $D\chi_F$ has no atoms), we write:
$$ P_\phi(F, \Omega) = \int_{\partial F} \phi(\nu(x)) d\mathcal{H}^1 = \sum_{i = 0}^{n} \phi(\nu_{\overline{q_i q_{i+1}}}) \mathcal{H}^1(\overline{q_i q_{i+1}}) $$
and
$$ P_\phi(F', \Omega) = \int_{\partial F'} \phi(\nu(x)) d\mathcal{H}^1 = \phi(\nu_{\overline{q_0 q_2}}) \mathcal{H}^1(\overline{q_0 q_2}) + \sum_{i = 2}^{n} \phi(\nu_{\overline{q_i q_{i+1}}}) \mathcal{H}^1(\overline{q_i q_{i+1}}).$$
As we assumed that $\alpha_i \nu_{\overline{q_i q_{i+1}}} \in I$ for every $i$, we see that also some positive multiple $\alpha$ of $\nu_{\overline{q_0 q_2}}$ belongs to $I$. Thus $(\nu_{\overline{q_0 q_2}} \cdotp \nu_{\overline{q_0 q_1}})^{-1} \alpha \nu_{q_0 q_1}$ and $(\nu_{q_0 q_2} \cdotp \nu_{\overline{q_1 q_2}})^{-1} \alpha \nu_{\overline{q_0 q_2}}$ belong to $I$. In this case, we have $\phi(\nu_{\overline{q_0 q_1}}) = (\nu_{\overline{q_0 q_2}} \cdotp \nu_{\overline{q_0 q_1}}) \phi(\nu_{\overline{q_0 q_2}})$. Similarly, $\phi(\nu_{\overline{q_1 q_2}}) = (\nu_{\overline{q_0 q_2}} \cdotp \nu_{\overline{q_1 q_2}}) \phi(\nu_{\overline{q_0 q_2}})$. Hence, we compare the two expressions and obtain

$$ P_\phi(F, \Omega) - P_\phi(F', \Omega) = \phi(\nu_{\overline{q_0 q_1}}) \mathcal{H}^1(\overline{q_0 q_1}) + \phi(\nu_{\overline{q_1 q_2}}) \mathcal{H}^1(\overline{q_1 q_2}) - \phi(\nu_{\overline{q_0 q_2}}) \mathcal{H}^1(\overline{q_0 q_2}) =  $$
$$ = \phi(\nu_{\overline{q_0 q_2}}) ((\nu_{\overline{q_0 q_2}} \cdotp \nu_{\overline{q_0 q_1}}) \mathcal{H}^1(\overline{q_0 q_1}) + (\nu_{\overline{q_0 q_2}} \cdotp \nu_{\overline{q_1 q_2}}) \mathcal{H}^1(\overline{q_1 q_2}) - \mathcal{H}^1(\overline{q_0 q_2})) = 0.$$
Thus the anisotropic perimeter of the sets $F$ and $F'$ bounded by the polygonal chains $q_0 q_1 q_2 ... q_{n+1}$ and $q_0 q_2 ... q_{n+1}$ respectively are the same, so inductively (with finitely many steps) we obtain that the anisotropic perimeter is the same as the anisotropic perimeter of $E$. \qed
\end{dd}

\begin{wn}\label{wn:nonuniquec1}
Let $F$ be a set such that its boundary is a piecewise $C^1$ curve from $p_1$ to $p_2$ such that the normal vector to $\partial F$ at each point lies in $N$. Then, as we can approximate it in the strict topology by sets whose boundaries are the polygonal chains as above, $F$ is a $\phi-$minimal set. \qed
\end{wn}

\begin{stw}\label{prop:construction}
Let $\Omega = B(0,1) \subset \mathbb{R}^2$ and suppose that $B_\phi(0,1)$ is not strictly convex. Then there exist boundary data $f \in C^{\infty}(\partial\Omega)$ such that the solution to Problem (\ref{problem}) is not unique.
\end{stw}

\begin{figure}[h]
    \includegraphics[scale=0.3]{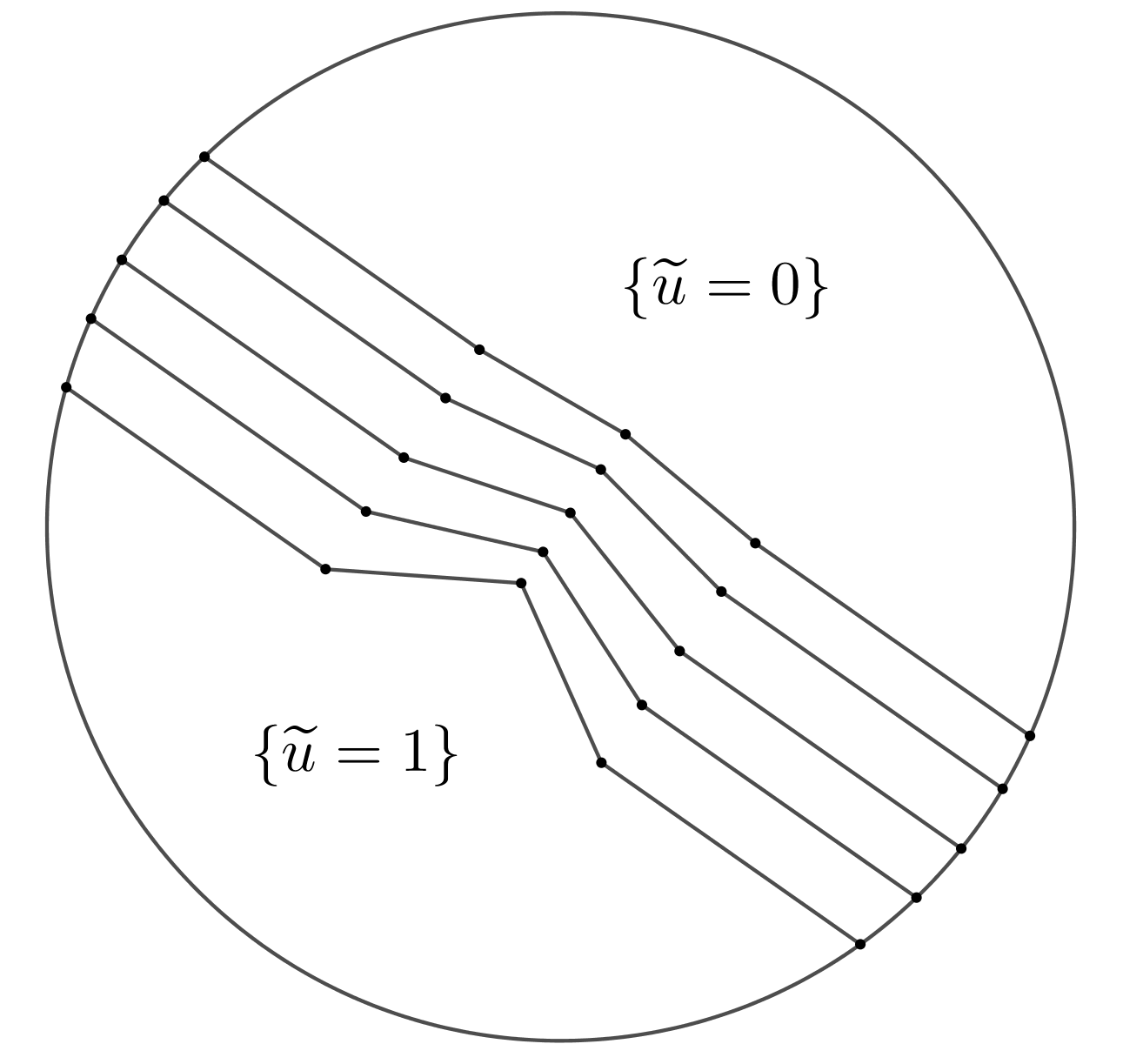}
    \caption{The construction of the competitor $\widetilde{u}$}
    \label{fig:funkcja}
\end{figure}

\begin{dd}
Let $I$ and $\nu_0$ be as above. Let $u \in C^{\infty}(\overline{\Omega})$ be such that it takes values in the interval $[0,1]$, all its level sets are line segments perpendicular to $\nu_0$ and $u$ is decreasing in the direction of $\nu_0$. Now, let $f = u|_{\partial \Omega}$. The preimage of every $t \in (0,1)$ consists of two points $p_1, p_2$ and the isotropic solution is such that the $t-$level set is the line segment $p_1 p_2$.

We want to perturb the line segment $p_1 p_2$ into a polygonal chain $\partial F$ consisting of four line segments, $\overline{p_1 q_1}, \overline{q_1 q_0}, \overline{q_0 q_2}, \overline{q_2 p_2}$ such that $q_1$ and $q_2$ lie on the line segment $\overline{p_1 p_2}$ and such that the vector $\nu_1$ normal to $q_1 q_0$ is sufficiently close to $\nu_0$ without changing the anisotropic perimeter. The length of the first and last of there line segments is $l_0$ and the length of the two in the middle is $l_1$: this way the vector $\nu_2$ normal to $\overline{q_0 q_2}$ is such that $\nu_1 + \nu_2$ is parallel to $\nu_0$, i.e. the triangle $\Delta q_0 q_1 q_2$ is an isosceles triangle. The sides of this triangle have lengths $l_1, l_1$ and $2 (\nu_0\cdotp \nu_1) l_1$. By Lemma \ref{lem:nonuniqueminimalsets} this polygonal chain has the same anisotropic perimeter as the original line segments, provided that all the normal vectors are such that a positive multiple of them lies in $I$. Hence, Proposition \ref{stw:anizobgg} implies that such $\widetilde{u}$ is another solution of Problem (\ref{problem}). This construction is summarized on Figure \ref{fig:funkcja}. \qed
\end{dd}

Now, we want to inspect the issue of regularity of the solutions to Problem (\ref{problem}) when $B_\phi(0,1)$ is not strictly convex. We present a series of examples, in which given a non-strictly convex metric integrand $\phi$ and an open bounded convex set $\Omega$ we construct a boundary datum $f \in C(\partial\Omega)$ such that there exists a solution that has regularity no better than $BV(\Omega) \cap L^\infty(\Omega)$.

\begin{stw}\label{stw:example}
Suppose that $\Omega$ is convex and $B_\phi(0,1)$ is not strictly convex. Then there exists a boundary datum, for which there exist solutions $u_1, u_2$ to Problem (\ref{problem}) such that $u_1 \notin W^{1,1}(\Omega)$ and $u_2 \notin SBV(\Omega)$.
\end{stw}

\begin{dd}
(1) Let $I \subset \partial B_{\phi}(0,1)$ be a line segment and let $\alpha \nu_0 \in \text{int } I$. Let $f \in C(\partial\Omega)$ be given by the formula

$$ f(x,y) = \nu_0 \, \cdotp \left( \! \!
  \begin{array}{cc}
    x \\
    y \\
  \end{array} \! \!
\right).$$
Then the function $u_0$ given by the same formula inside $\Omega$ has the prescribed trace $f$ and all its superlevel sets are line segments; by Proposition \ref{prop:linesegment} they are $\phi-$minimal. Thus by Proposition \ref{stw:anizobgg} $u_0$ is a function of $\phi-$least gradient.

We are going to modify the function $u_0$ on a compact subset of $\Omega$, so that the trace remains unchanged. We will modify carefully the level sets of $u_0$. To this end, we are going to use Proposition \ref{lem:nonuniqueminimalsets} to construct as superlevel sets of $u_1$ another sets with the same perimeter and trace.

Take any line $l$ perpendicular to $\nu_0$ intersecting $\Omega$. Choose two points $q_1, q_2 \in l$. Then, take two lines $m, m'$ perpendicular to $l$, intersecting $l$ inside $\Omega$, such that the points $q_1, q_2$ lie between $m$ and $m'$. Now, we choose two other lines $l', l''$ parallel to $l$ such that $l$ lies between $l'$ and $l''$ and such that they are sufficiently close to $l$, so the following condition is fulfilled: let $p_1 = l' \cap m, p_2 = l' \cap m', p_3 = l'' \cap m, p_4 = l'' \cap m'$. Then a positive multiple of the normal vector to each of the line segments $\overline{p_i q_j}$ belongs to $I$. 

We construct each of the level sets of $u_1$ in the following way: outside the region bounded by $l$ and $l'$ we define it to equal $u_0$. Now, fix any line $k$ parallel to $l$ lying between $l'$ and $l''$. Then the level set is a polygonal chain along this line until the point $k \cap m$, then it is a line segment from this point to $q_1$, then the line segment $\overline{q_1 q_2}$, then again from $q_0$ to $k \cap m'$ and then again along $k$. This construction is shown on Figure \ref{fig:nieW11}.

\begin{figure}[h]
    \includegraphics[scale=0.3]{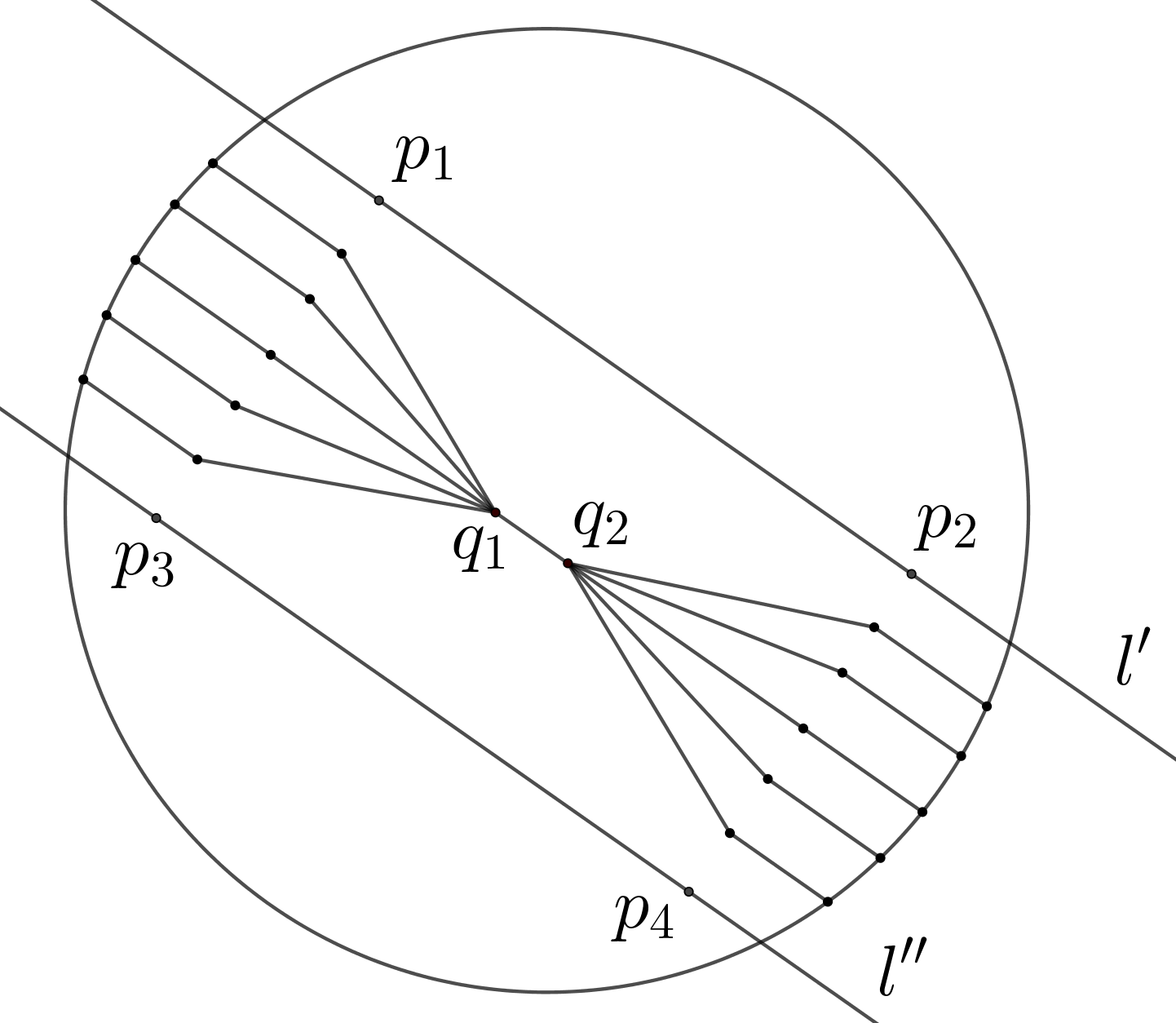}
    \caption{The construction of another minimizer $u_1$}
    \label{fig:nieW11}
\end{figure}
By Lemma \ref{lem:nonuniqueminimalsets} the $\phi-$perimeter of each of the level sets of $u_1$ is the same as the $\phi-$perimeter of each of the level sets of $u_0$, as the normal vector to each of the line segments lies in $I$. By Proposition \ref{stw:anizobgg} $u_1$ is a function of $\phi-$least gradient. We notice that the constructed minimizer has a constant jump along the line segment $\overline{q_1 q_2}$; therefore it fails the ACL characterization and does not belong to $W^{1,1}(\Omega)$.

\vspace{1mm}

(2) Again, we are going to modify $u_0$ on a compact subset of $\Omega$ and the trace remains unchanged. Take any line $l$ perpendicular to $\nu_0$ intersecting $\Omega$. Take four lines $m_1, m_2, m_3, m_4$ perpendicular to $l$ intersecting $l$ inside $\Omega$ (in this order along $l$). Now, take a parallel line $l'$ sufficiently close to $l$, so that the following condition is satisfied: let $\{ p_i \}$ be all eight possible intersections between the lines $l$ or $l'$ and the lines $m_k$. Then a positive multiple of the normal vector to the each of the line segments $\overline{p_i p_j}$, where $\overline{p_i p_j}$ does not lie on $m_k$ for some $k$, belongs to $I$.

We construct each of the level sets of $u_2$ in the following way: without loss of generality, we may assume that $u_0 \equiv 0$ on $l$ and $u_0 \equiv 1$ on $l'$. Let $g$ be the Cantor stairs function. The level sets of $u_2$ for $t \notin (0,1)$ are the same as for $u_0$. Now, fix $t \in (0,1)$. Let $l_t$ denote the line parallel to $l$ corresponding to value $t$ of the minimizer $u_0$. Then the level set $\{ u = t \}$ is as follows: firstly along $l_t$ up to the intersection with $m_1$; then the line segment $[l_t \cap m_1, l_{\frac{t+g(t)}{2}} \cap m_2]$; then along $l_{\frac{t+g(t)}{2}}$ up to the intersection with $m_3$; then the line segment $[l_{\frac{t+g(t)}{2}} \cap m_3, l_t \cap m_4]$; finally, again along $l_t$. This construction is shown on Figure \ref{fig:nieSBV}.

Again, by Lemma \ref{lem:nonuniqueminimalsets} and Proposition \ref{stw:anizobgg} $u_2$ is a function of $\phi-$least gradient. We notice that the derivative of the constructed minimizer in the rectangle bounded by the lines $l, l'', m_2, m_3$ is a continuous measure which is not absolutely continuous; thus $u_2 \notin SBV(\Omega)$.

\begin{figure}[h]
    \includegraphics[scale=0.3]{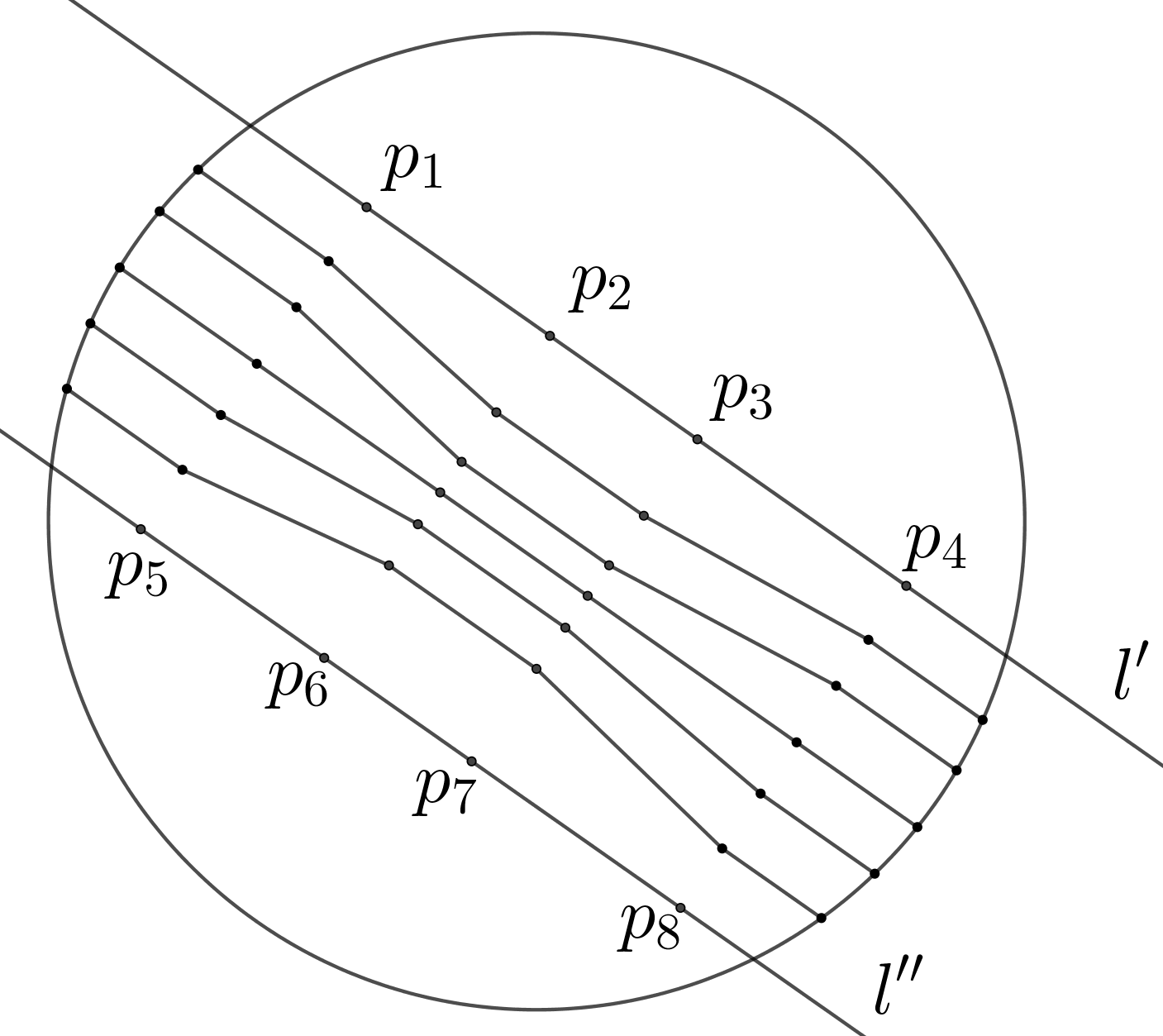}
    \caption{The construction of another minimizer $u_2$}
    \label{fig:nieSBV}
\end{figure}

\qed
\end{dd}

Finally, we turn our attention to the barrier condition. The next Proposition justifies the need for Theorems \ref{tw:existence} and \ref{tw:existence2}, as we cannot apply the existing theory from \cite{JMN} to obtain existence of minimizers to Problem (\ref{problem}).

\begin{stw}\label{stw:barriernothold}
Suppose that $\Omega$ has $C^1$ boundary and $B_\phi(0,1)$ is not strictly convex. Then $\Omega$ does not satisfy the barrier condition.
\end{stw}

\begin{dd}
By Proposition \ref{prop:linesegment} a line segment is a $\phi-$minimal surface. Therefore if $\Omega$ is not strictly convex, it does not satisfy the barrier condition.

Now, take $\Omega$ to be a strictly convex set. As $B_\phi(0,1)$ is not strictly convex, there exists a line segment $I \subset \partial B_\phi(0,1)$. Fix $\nu_0$ such that $\alpha \nu_0 \in $ int $I$. Take $x_0 \in \partial\Omega$ such that the normal vector to $\partial\Omega$ at $x_0$ has direction $\nu_0$. Take $\varepsilon$ small enough, so that a positive multiple of the normal vector $\nu$ at every $x \in \partial\Omega \cap B(x_0, \varepsilon)$ lies in $I$.

Provided $\varepsilon$ is small enough, the set $\partial\Omega \cap \partial B(x_0,\varepsilon)$ consists of two points $x_1, x_2$. Take $V$ to be the open connected set bounded by $\partial\Omega$ and the two line segments $[x_1, x_0]$ and $[x_2, x_0]$. Then $V \backslash B(x_0, \varepsilon) = \Omega \backslash B(x_0, \varepsilon)$ and by Lemma \ref{lem:nonuniqueminimalsets} $V$ is a minimal set. However, $\partial V \cap \partial \Omega \cap \partial B(x_0, \varepsilon) = \{ x_0 \} \neq \emptyset$, so the barrier condition is not satisfied. \qed
\end{dd}

However, if $\partial\Omega$ has corners, it is still possible that the barrier condition is satisfied. This depends on the balance between the width of an angle corresponding to a flat part of $\partial B_\phi(0,1)$ and the width of the corner. For instance, it is easy to see that for $l^1$ anisotropy no open set satisfies the barrier condition; the next Example shows that if the flat part of $\partial B_\phi(0,1)$ is small enough, the barrier condition may hold for properly chosen $\Omega$.

\begin{prz}
Let $\phi$ be an anisotropic norm such that $\partial B_\phi(0,1)$ has only two flat parts on the boundary corresponding to angles $(\frac{1\pi}{8}, \frac{3\pi}{8})$ and $(\frac{9\pi}{8}, \frac{11\pi}{8})$ in the polar coordinates on the plane. Let $\Omega$ be an open strictly convex set, symmetric with respect to the line $y = x$, such that it is $C^\infty$ except for two corners at $\pm (1,1)$, such that the angle of incidence of $\partial\Omega$ at $\pm (1,1)$ to the line $y = x$ is $\frac{\pi}{16}$ and the angle of incidence of $\partial \Omega$ to any line of the form $y = x + a$ is bounded from above by $\frac{\pi}{16}$. Then the barrier condition is satisfied, as any $\phi-$minimal surface with an angle of incidence to $y=x+a$ smaller than $\frac{\pi}{8}$ is a line segment. 

The proof of Proposition \ref{stw:barriernothold} fails, as we cannot take a $x \in \partial\Omega$ such that in its neighbourhood the normal vector to $\partial\Omega$ has direction corresponding to a flat part on the boundary; here, there are only two isolated points, $\pm (1,1)$, with such normal vectors.
\end{prz}

As the solutions to Problem (\ref{problem}) may be not unique if $B_\phi(0,1)$ is not strictly convex, we may ask if the solution with minimal $L^1$ norm exists and has any additional regularity (this issue has been discussed for the isotropic case with discontinuous boundary data in \cite{Gor2}). However, using a technique similar as in the previous Proposition, we may prove that existence of such solutions fails for sets with $C^1$ boundary.

\begin{stw}
Suppose that $\Omega \subset \mathbb{R}^2$ has $C^1$ boundary. Suppose that $B_\phi(0,1)$ is not strictly convex. Then there exist $C^\infty$ boundary data such that there is no minimizer of Problem (1) with minimal $L^1$ norm.
\end{stw}

\begin{dd}
As before, we may assume that $\Omega$ is strictly convex and $\nu_0, x_0, \varepsilon, x_1$ and $x_2$ are as in the proof of the previous Proposition. Take $f \in C^\infty(\partial\Omega)$ satisfying the following conditions: $f = 0$ on $\partial\Omega \backslash B(x_0, \varepsilon)$; $f(x_0) = 1$; $f$ is strictly monotone on the arcs $(x_1, x_0)$ and $(x_0, x_2)$ on $\partial\Omega$; $f$ is one-dimensional in the direction of $\nu_0$, i.e. $f(x,y) = \widetilde{f}(\nu_0 \cdotp (x,y))$.

Define the functions $u_n$ in the following way: let $y_n \in (x_0 + \nu_0 \mathbb{R}) \cap \Omega \cap B(x, \varepsilon)$ and $y_n \rightarrow x_0$. For $y \in [y_n, x_0]$ we fix $u_n(y) = t$, where $y = (1 - t)y_n + tx_0$. Denote by $x_1^t$ and $x_2^t$ the two elements of $f^{-1}(t)$ for $t \in (0,1)$. Then, we fix $u_n = t$ on a shifted boundary of $\Omega$, namely $\partial\Omega + (y_n - x_0)$ in a smaller ball $B(x_0, \varepsilon')$ and $u_n = t$ on the line segments $[x_1^t, y_1^t], [y_2^t, x_2^t]$, where $y_1^t$ and $y_2^t$ are the two points of $(\partial_\Omega + (y_n - x_0)) \cap \partial B(x_0, \varepsilon')$. This sequence converges to $0$ everywhere in $\Omega$ (as $y_n \rightarrow x_0$) and by Lemma \ref{lem:nonuniqueminimalsets} each $u_n$ is a function of least gradient. By construction also $Tu_n = f$. However, $\| u_n \|_{L^1(\Omega)} \rightarrow 0$ and $u \equiv 0$ is not a minimizer of Problem (\ref{problem}), as it does not satisfy the boundary condition. Thus the minimum $L^1$ norm among minimizers is not attained. \qed
\end{dd}

\begin{figure}[h]
    \includegraphics[scale=0.3]{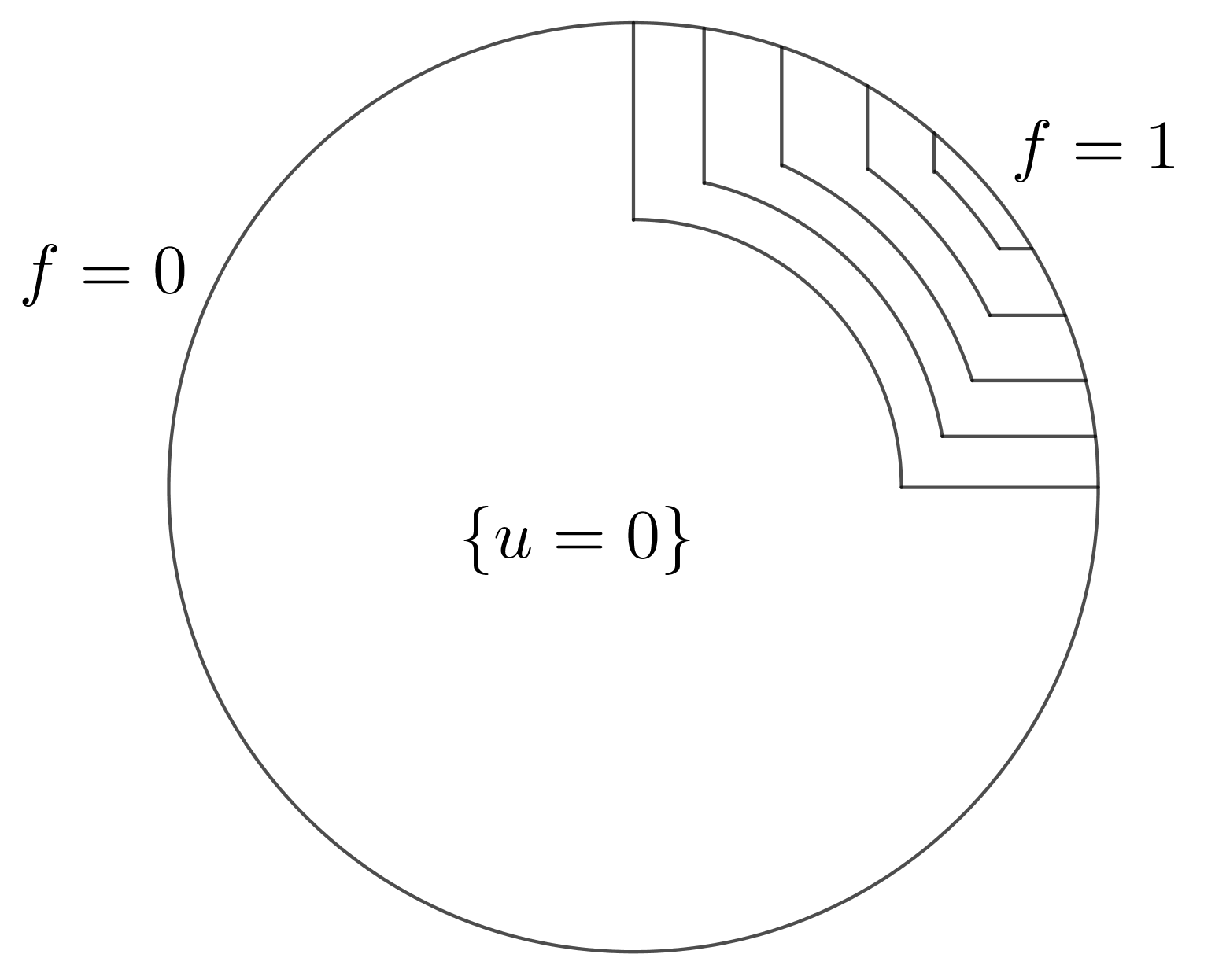}
    \caption{Minimizers with arbitrarily small $l^1$ norm}
    \label{fig:malanorma}
\end{figure}

The above construction was given to show that for every $\phi$ such that $B_\phi(0,1)$ is not strictly convex there exist nonzero boundary data such that the infimum of the $L^1$ norms of the solutions equals $0$. An example illustrating this phenomenon for the $l_1$ anisotropy is shown on Figure \ref{fig:malanorma}.

\section{Conclusion}

In this final Section we want to discuss what implications do the results from previous Sections have for the validity of several results from the isotropic case also in the anisotropic case. We will recall these results and shortly discuss the methods used in their proofs to see that they also hold in the anisotropic case with properly modified proofs.

The first result concerns the existence of minimizers of Problem (\ref{problem}) if the boundary data lie in the class $BV(\partial\Omega)$. For the isotropic version, see \cite[Theorem 1.1]{Gor1}.

\begin{tw}\label{tw:gor1}
Let $\Omega$ be an open bounded strictly convex set with $C^1$ boundary. Suppose that $B_\phi(0,1)$ is strictly convex. Then for $f \in BV(\partial\Omega)$ there exists a solution to Problem (\ref{problem}).
\end{tw}

\begin{dd}
The proof of Theorem \ref{tw:gor1} in the isotropic case has the following outline: we approximate the boundary data $f \in BV(\partial\Omega)$ with smooth functions $f_n$ in the strict topology (in the isotropic norm). We take advantage of the fact that on the one-dimensional manifold, such as $\partial\Omega$, the (isotropic) perimeter is a natural number (or is infinite), so $\partial \{ f \geq t \}$ consists of finitely many points for almost all $t \in \mathbb{R}$. We construct the sets $\{ u_n \geq t \}$; their boundaries, in the anisotropic case by \ref{stw:regularityofminimalsets}, are unions of (finitely many) line segments connecting points from $\partial \{ f \geq t \}$. We construct the minimizers $u_n$ as in the proof of Theorem \ref{tw:gor1} and their convergence to a minimizer is guaranteed by Theorem \ref{tw:gammaconvergence}.  Thus, due to Proposition \ref{stw:regularityofminimalsets}, the Theorem is adaptable to the anisotropic case with essentially the same proof. \qed
\end{dd}

The second (\cite[Theorem 1.2]{Gor1}) and third (\cite[Theorem 1.1]{Gor2}) result concern the structure of minimizers for arbitrary boundary data. The third result can also be understood as a uniqueness-type result; for discontinuous boundary data, the minimizers need not be unique, see \cite[Example 2.7]{Maz}, but this result shows that the minimizers may differ only on level sets of positive Lebesgue measure. Obviously, Example \ref{stw:example} shows that these results do not hold in $B_\phi(0,1)$ is not strictly convex.

\begin{tw}
Let $\Omega$ be an open bounded convex set. Suppose that $B_\phi(0,1)$ is strictly convex. Let $u$ be a minimizer of Problem (\ref{problem}). Then $u_j = u_c + u_j$, where $u_c$ is continuous, $u_j$ has only jump-type derivative and this decomposition is unique up to an additive constant.
\end{tw}

\begin{tw}
Let $\Omega$ be an open bounded convex set. Suppose that $B_\phi(0,1)$ is strictly convex. Let $u,v$ be two minimizers of Problem (\ref{problem}). Then $u - v$ is a locally constant function.
\end{tw}

The proof of these results involves mostly heavy regularity theory for minimal sets, a (strong) maximum principle for minimal surfaces (inside $\Omega$) and a weak maximum principle, i.e. the fact that two connected components of $\partial \{ u > t \}$ cannot intersect on the boundary of $\Omega$; as the regularity theory is irrelevant due to Proposition \ref{stw:strict} and the last result is adaptable to the anisotropic case, these results hold also in the anisotropic case. 

Finally, it is worth noting that in the isotropic case, if $\Omega$ is only convex, we may still obtain existence of minimizers, if $f$ satisfies some additional admissibility conditions, see \cite{RS}; in principle, due to results from Section 4, the methods used there are also adaptable in the anisotropic case if $B_\phi(0,1)$ is strictly convex. However, justifying this for suitable anisotropic version of the admissibility conditions goes beyond the scope of this paper.

\begin{ak}
I would like to thank my PhD advisor, Piotr Rybka, for fruitful discussions about this paper. This work was partly supported by the research project no. 2017/27/N/ST1/02418, "Anisotropic least gradient problem", funded by the National Science Centre.
\end{ak}

\bibliographystyle{plain}
\bibliography{WG-references}

\end{document}